\documentclass[10pt]{article}
\usepackage{amsmath,amsthm,amscd,amssymb}
\newtheorem{df}{Definition}[section]
\newtheorem{thm}{Theorem}[section]
\newtheorem{prop}{Proposition}[section]
\newtheorem{lm}{Lemma}[section]

\newtheorem{remark}{Remark}[section]
\newtheorem{fact}{Fact}[section]
\newtheorem{cor}{Corollary}[section]

\title{The Taylor expansion of Ruelle L-function at the origin and the
Borel regulator}
\author{Ken-ichi SUGIYAMA\thanks{Cooresponding address : Department of
Mathematics and Informatics, Faculty of Science, Chiba University, 1-33
Yayoi-cho Inage-ku, Chiba 263-8522, Japan. e-mail address : sugiyama@math.s.chiba-u.ac.jp}\\
Department of Mathematics and Informatics, Faculty of Science\\
Chiba University, Japan.}

\begin{document}
\maketitle
\begin{abstract}
We will prove that Ruelle L-function for a cuspidal unitary local system
 on an odd dimensional complete hyperbolic manifold with finite volume satisfies a functional equation and
 an analogue of Riemann hypothesis.  We will also compute its Laurent expansion at the
 origin and will prove that the second coefficient coincides with a rational multiple
 of the volume up to a certain contribution from cusps. Moreover if the dimension is
 three we will identify the leading coefficient. Both of them will be interpreted as a period of a certain element of K-group of 
${\mathbb C}$. Also a relation with the $L^2$-torsion will be discussed.
\par\vspace{5pt}
\end{abstract}
\section{Introduction}
Researches of an L-functions may be roughly classfied in the
following three subjects:

\begin{enumerate}
\item a functional equation,
\item (Riemann hypothesis) a distribution of zeros and poles,
\item an arithmetic or a geometric meaning of special values. 
\end{enumerate}

For example let us consider the zeta function for a number field $F$:
\[
 \zeta_F(s)=\prod_{\frak{P}}(1-e^{-s\log N(\frak{P})})^{-1}.
\]
Here $\frak{P}$ runs thorough all prime ideals and $N(\frak{P})$ is the
norm. Then $\zeta_F(s)$ satisfies a functional
equation and the Riemann hypothesis is still a far reaching
problem. It has a zero at $s=0$ of order $r_1+r_2-1$ where
$r_1$ (resp. $r_2$) is the number of real (resp. complex) places and the
leading coefficient of the Taylor expansion is given by
\[
 \lim_{s\to 0}s^{-(r_1+r_2-1)}\zeta_F(s)=-\frac{\sharp{\rm Pic}({\mathcal
 O}_F)}{\sharp({\mathcal O}^{\times}_F)_{tors}}\cdot R.
\]
$R$ is the covolume of the image of ${\mathcal
O}^{\times}_F\oplus{\mathbb Z}$ by the classical regulator which is
defined by the logarithmic function:
\[
 {\mathcal
O}^{\times}_F\oplus{\mathbb Z} \stackrel{r_{1,F}}\to {\mathbb R}^{r_1+r_2}.
\]
The
observation that $K_1({\mathcal O}_F)$ is isomorphic to ${\mathcal
O}_F^{\times}$ and the fact that the order of $\zeta_F(s)$ at $s=1-l$ is
equal to the dimension $d_l$ of $K_{2l-1}(F)\otimes{\mathbb Q}$ for
$l\geq 2$ lead
Lichtenbaum to a conjecture; There should be a map
\[
 K_{2l-1}(F)\stackrel{r_{l,F}}\to {\mathbb R}^{d_l}.
\]
such that
\[
 \lim_{s\to 1-l}(s+1-l)^{-d_l}\zeta_F(s)=vol({\mathbb R}^{d_l}/r_{l,F}(K_{2l-1}(F))).
\]

This was solved by Borel. He has also constructed a map
\[
 K_{2l-1}({\mathbb C})\stackrel{r_l}\to {\mathbb R},
\] 
and an each of $r_{l,F}$ and $r_l$ is referred as {\it the Borel
regulator}. In this paper, under a certain condition, we will show that Ruelle L-function
for a unitary local system on an odd dimensional hyperbolic manifold
(especially a threefold) with finite volume carries similar
properties.\\

Let $X$ be a complete hyperbolic $d$-fold ($d=2n+1$, $n\geq 1$) of
finite volume. Thus it is a quotient of the Poincar\'{e} upper half space
${\mathbb H}^d$ by a torsion free discrete subgroup $\Gamma$ in ${\rm
SO}_0(d,1)$, the connected component of the isometry groups ${\rm
SO}(d,1)$ of ${\mathbb H}^d$. Notice that there is the natural bijection between a
set of hyperbolic conjugacy classes $\Gamma_{hyp}$ of $\Gamma$ and a set of closed
geodesics of $X$. By this the length $l(\gamma)$ of a hyperbolic
conjugacy class $\gamma$ is defined to be one of the corresponding
closed geodesic. A closed geodesic will be referred as {\it prime} if it is not
a positive multiple of an another one. Using the bijection we define a
subset $\Gamma_{prim}$ of hyperbolic conjugacy classes which
consists of elements corresponding to prime closed geodesics. Let $\rho$
be a unitary representaion of $\Gamma$ with degree $r$,  i.e. the dimension
of the representation space $V_{\rho}$ is $r$.  
Now Ruelle L-function is defined to be
\[
 R_{X}(z,\rho)=\prod_{\gamma\in \Gamma_{prim}}\det[1-\rho(\gamma)e^{-zl(\gamma)}]^{-1}.
\]
It absolutely converges if ${\rm Re}s>2n$ and is
meromorphically continued to the whole plane (\cite{Park2007}, see also
\S 2.4). Hereafter, otherwise mentioned, we will assume that 
${\rho}$ is {\it cuspidal} (see \S 2). \\

We will show that $R_{X}(z,\rho)$ satisfies a functional equation 
\[
 R_{X}(z,\rho)\cdot R_{X}(-z,\rho)^{-1}=\exp[\frac{vol(X)}{\pi}Y(z)+4\sum_{j=0}^n(-1)^j\delta(X,\rho)z],
\]
where $Y(z)$ is a polynomial of rational coefficients which vanishes at
$z=0$ and $\delta(X,\rho)$ is a certain constant determined by special
values of Epstein L-functions of the fundamental groups at
cusps (see \S 2.3). Notice that here and hereafter if $X$ is closed, since it has no cusp,
$\delta(X,\rho)$ does not appear. It will be also shown that its zeros and poles are located on 
\[
 \{z\in {\mathbb C}\,|\, {\rm Re}z=-n,-(n-1),\cdots,n-1,n\},
\]
except for finitely many of them.
For example if $d=3$, i.e. $n=1$, 
\[
 R_{X}(z,\rho)\cdot R_{X}(-z,\rho)^{-1}=\exp[\frac{2r}{\pi}vol(X)(\frac{z^3}{3}-3z)].
\]
Its logarithmic derivative 
\[
 r_X(z,\rho)+r_X(-z,\rho)=\frac{2r}{\pi}vol(X)(z^2-3), \quad
 r_X(z,\rho)=\frac{d}{dz}\log R_{X}(z,\rho), 
\]
may be compared to the functional equation of Weil conjecture. In fact Hasse-Weil's congruent zeta function of a smooth projective variety
$M$ with dimension $m$ over a finite field ${\mathbb F}_q$ is defined to be
\[
 \zeta_M(z)=\prod_{P\in |M|}(1-q^{-z{\rm deg}(x)})^{-1},
\]
where $|M|$ is the set of closed points of $M$ and ${\rm deg}(x)$ is the
extension degree of the resdue field of $x$ over ${\mathbb F}_q$. Weil
conjecture implies that
\[
 \sigma_M(z)=\frac{d}{dz}\log_q \zeta_M(z+\frac{m}{2}),
\]
satisfies
\[
 \sigma_M(z)+\sigma_M(-z)=\chi(M),
\]
where $\chi(M)$ is the Euler characterisitc of $M$. Thus replacing $\chi(M)$ by $2rvol(X)(z^2-3)/\pi$ (which may be
not so absurd if we think about Gauss-Bonnet's formula), we find that
the logarithmic derivative of Ruelle L-function for a cuspidal unitary local system on a hyperbolic
threefold satisfies a functional equation similar to $\sigma_M(s)$.\\

Let us expand $R_X(z,\rho)$ at the origin which is a symmetric point of
the functional equation:
\[
 R_X(z,\rho)=c_0z^h(1+c_1z+\cdots),\quad c_0\neq 0.
\]
We are interested in the coefficients $c_0$ and $c_1$.
\begin{thm} $c_1-2\sum_{j=0}^n(-1)^j\delta(X,\rho)$ is a rational multiple of $vol(X)/{\pi}$.
\end{thm}
For example if $d=3$, we will show 
\[
 c_1=-3r\frac{vol(X)}{\pi}.
\]
Combining with the results of Goncharov(\cite{Goncharov}) this yields
\begin{cor}
There is an element $\gamma_X\in K_{2n+1}({\mathbb C})$ called {\it the
 Borel element} so that $c_1-2\sum_{j=0}^n(-1)^j\delta(X,\rho)$ is a rational
 multiple of $r_{n+1}(\gamma_X)/{\pi}$.  
\end{cor}
In \S4 we will recall a construction of
$\gamma_X$ for a closed hyperbolic threefold and will show 
\[
 c_1=-\frac{3r}{\pi}\cdot r_2(\gamma_X).
\]
It is natural to expect that the leading coefficient $c_0$ has a similar
interpretation. In fact it is true at least if $d=3$.
\begin{thm} Let $h^p(X,\rho)$ be the dimension of $H^p(X,\rho)$. Suppose
 that $d=3$. Then $R_X(z,\rho)$ has a zero at the origin of order
 $2h^1(X,\rho)$ and 
\[
 c_0=(\tau^{*}(X,\rho)\cdot {\rm Per}(X,\rho))^2.
\]
\end{thm}
Here $\tau^{*}(X,\rho)$ and ${\rm Per}(X,\rho)$ are {\it the modified
Franz-Reidemeister torsion} and {\it the period} of $(X,\rho)$,
respectively. (See \S 3.4.) (If $h^1(X,\rho)$ is zero $\tau^{*}(X,\rho)$
is the usual Franz-Reidemeister torsion $\tau(X,\rho)$.) If $X$ is
closed the theorem has been already proved by Fried (\cite{Fried}). In
fact he has proved it for a closed odd dimensional hyperbolic
manifold. \\

Suppose $d=3$ and that 
$h^1(X,\rho)$ vanishes. Let us fix a triangulation of $X$ and a unitary basis ${\bf
e}=\{{\bf e}_1,\cdots,{\bf e}_r\}$ of $V_{\rho}$. By Poincar\'{e} duality we know that all $H^p(X,\rho)$
vanishes for all $p$ and therefore the cochain
complex $C^{\cdot}(X,\rho)$ is acyclic. Then Milnor has constructed an
element $\tau(X,\rho,{\bf e})$ in $K_1({\mathbb C})\simeq {\mathbb
C}^{\times}$ which is referred as {\it the Milnor element} and has shown
(\cite{MilnorI}):
\[
 \log \tau(X,\rho)=2\pi r_1(\tau(X,\rho,{\bf e})).
\]
Summarizing there are rational numbers $\alpha$ and $\beta$ such that
\[
 \log R_X(0,\rho)=\alpha\pi r_1(\tau(X,\rho,{\bf e})), \quad \frac{d}{dz}\log R_X(z,\rho)|_{z=0}=\frac{\beta}{\pi}\cdot r_2(\gamma_X).
\]
Thus replacing the logarithmic derivative by a shift:
\[
 f=f(z) \to f^{[k]}(z)=f(z-k),
\]
our formula will correspond to
Lichtenbaum conjecture. \\

In \cite{Park2007}, Park has obtained
\begin{fact}
Let $X$ be an odd dimensional complete hyperbolic manifold with finite
 volume and $\rho$ a unitary local system on $X$ which may not satisfy
 the cuspidal condition. Then the leading coefficient of the Laurent
 expansion at the origin is $\exp(-\zeta^{\prime}_X(0,\rho))$ where
 $\zeta_X(s,\rho)$ is the spectral zeta function (see \S 3.4).
\end{fact}
By Hodge theory $H^p(X,\rho)$ is isomorphic to ${\rm Ker}\Delta^p_X$, the kernel of Hodge Laplacian
$\Delta^p_X$. They are subspaces of $C^{\cdot}(X,\rho)$ and $L^2(X,\Omega^p(\rho))$, the space of square integrable sections of $p$-forms twisted by $\rho$
on $X$, respectively. (Here notice that $H^p(X,\rho)$ is isomorphic to
the kernel of combinatric Laplacian acting on $C^{p}(X,\rho)$, see \S3.4
.) Using this two metrics will be defined on the determinant
line bundle $\det H^{\cdot}(X,\rho)$. One is {\it Franz-Reidemeister metric}
which is defined in terms of the combinatric $L^2$-norm
$|\cdot|_{l^2,X}$ induced from the natural metric on $C^{\cdot}(X,\rho)$ and the
modified Franz-Reidemeister torsion: 
\[
 ||\cdot||_{FR}=|\cdot|_{l^2,X}\cdot \tau^{*}(X,\rho)^{1/2}.
\]
The other is {\it Ray-Singer metric}: 
\[
 ||\cdot||_{RS}=|\cdot|_{L^2,X}\cdot\exp(-\frac{1}{2}\zeta^{\prime}_{X}(0,\rho)),
\]
where $|\cdot|_{L^2,X}$ is the analytic $L^2$-norm derived from the inner
product on $L^2(X,\Omega^p(\rho))$. Since by definition ${\rm
Per}(X,\rho)$ is $|\cdot|_{l^2,X}/|\cdot|_{L^2,X}$, {\bf
Theorem 2.1} is reduced to show the following Cheeger-M\"{u}ller type theorem.
\begin{thm}
\[
 ||\cdot||_{FR}=||\cdot||_{RS}.
\]
\end{thm}
For a convenience we will give a proof of {\bf Fact 1.1} in {\bf
Appendix} under an assumption
that $d=3$ and $\rho$ is cuspidal. In particular the last assumption implies that
it is not neccessary to take care of the scattering term in
Selberg trace formula and one can prove the desired result just following
Fried's argument (\cite{Fried}). 

{\bf Acknowledgment.} The author express heartly gratitude to Professor
J. Park who kindly show him a preprint \cite{Park2007} which is
indispensable to finish this work.

\section{The second coefficient}
Let \[
 {\mathbb H}^d=\{(x_1,\cdots,x_{d+1})\,|\, x_1^2+\cdots+x_d^2-x_{d+1}^2=-1,\,x_{d+1}>0\}
\]
be the hyperbolic space form $(d=2n+1,\, n\geq 1)$.  We will choose its
origin to be ${\bf
o}=(0,\cdots,0,1)$ . The connected component $G={\rm SO}_o(d,1)$ of its
isometry group ${\rm SO}(d,1)$ transitively acts on ${\mathbb H}^d$ and the isotropy
subgroup at ${\bf o}$ is a maximal compact
subgroup $K={\rm SO}(d)$. Thus we have a surjective map
\[
 G\stackrel{\pi}\to {\mathbb H}^d, \quad \pi(g)=g\cdot{\bf o},
\]
which induces a diffeomorphism 

\begin{equation}
G/K\simeq {\mathbb H}^d.
\end{equation}
Let ${\frak g}$ be the Lie algebra of $G$ and $\theta$ the Cartan
involution. We define {\it the normalized Cartan-Killng form} to be
\[
 (X,Y)=-\frac{1}{4\pi}{\rm Tr}(ad X\circ ad (\theta Y)), \quad X,Y\in
 {\frak g}.
\]
The Cartan involution provides a decomposition of the Lie algebra
${\frak g}$ of $G$: 
\[
 {\frak g}={\frak k}\oplus{\frak h},
\]
where ${\frak k}$ and ${\frak h}$ are the $+1$- and $-1$-eigenspaces,
respectively. ${\frak h}$ may be
identified with the tangent space of ${\mathbb H}^n$ at the origin and the normalized Cartan-Killng form defines a Riemannian metric on
${\mathbb H}^n$ with constant curvature $-1$. Let ${\frak a}$ be the
maximal abelian subalgebra of ${\frak h}$ and
$\beta$ the positive restricted root of $({\frak g},{\frak a})$. Then
${\frak a}$ is one dimensional and let us
choose $H\in {\frak a}$ satisfying $\beta(H)=1$. Then the Lie subgroup $A$
of ${\frak a}$ is isomorphic to ${\mathbb R}$ by a
map:
\begin{equation}
 {\mathbb R} \simeq A,\quad t \mapsto \exp(tH).
\end{equation}
Let ${\frak n}$ be the positive root space of $\beta$ and $N=\exp({\frak
n})$ the associated Lie subgroup. Then using the Iwasawa decomposition
\[
 G=KAN,
\]
we introduce a Haar measure on $G$ by
\[
 dg=a^{2\rho}dk\cdot da\cdot dn.
\]
Here $\rho=n\beta$ is the half sum of potitive roots of $({\frak
g},{\frak a})$ and $a^{2\rho}=\exp(2\rho(\log a))$. $dk$ is the Haar
measure on $K$ whose total mass is one and $da$ is the push forward of
the Lebesgue measure on ${\mathbb R}$ by (2). The
volume form $dn$ of
$N$ is induced by the normalized Cartan-Killng form. Let $M\simeq {\rm SO}(d-1)$ be the centralizer of $A$ in $K$ and $P=MAN$
a proper parabolic subgroup.\\

 Let $X$ be a complete hyperbolic $d$-fold with finite volume and
$\{\infty_1,\cdots,\infty_{h}\}$ be the cusps. Thus it is a quotient of
${\mathbb H}^d$ by a torsion free discrete subgroup $\Gamma$ in $G$. A conjugate $P_{\nu}=g_{\nu}Pg^{-1}_{\nu}$ $(g_{\nu}\in G)$ corresponds
to a cusp $\infty_{\nu}$ and {\it the fundamental group} $\Gamma_{\nu}$
at $\infty_{\nu}$ is defined to be
\[
 \Gamma_{\nu}=\Gamma\cap P_{\nu}.
\]
We will normalize so that $\infty_1=\infty$ and $g_1$ is the
identity. 
Since $\Gamma$ is torsion free $\Gamma_{\nu}$ is equal to $\Gamma \cap N_{\nu}$
($N_{\nu}=g_{\nu}Ng_{\nu}^{-1}$) which is a lattice in ${\mathbb R}^{2n}$.
Let $\rho$ be a unitary representation of $\Gamma$ of degree $r$ and
$\rho_{\nu}$ its restriction to $\Gamma_{\nu}$. Since $\Gamma_{\nu}$ is
abelian $\rho_{\nu}$ is decomposed into
a direct sum of characters:
\[
 \rho_{\nu}=\oplus_{i=1}^{r}\chi_{\nu,i}.
\]
Throughout the paper we will assume that ${\rho}$ is {\it cuspidal},
i.e. {\it none} of
$\{\chi_{\nu,i}\}_{1\leq i \leq r, 1\leq \nu \leq h}$ is trivial. This
terminology will be justified in {\bf Lemma 3.1}.\\

Let $\Omega^j_X$ be the vector bundle of $j$-forms on $X$ and
$\Omega^j_X(\rho)$ its twist by $\rho$. Then the pullback $\Omega^j$ of
$\Omega^j_X$ on ${\mathbb H}^d$ is a homogeneous vector bundle. In fact
let $\xi$ be the standard action of ${\rm SO}(d)$ on ${\mathbb R}^d$ and ${\rm
SO}(d)\stackrel{\xi_j}\to{\rm GL}(\wedge^j{\mathbb R}^d)$ its exterior
product.  Then $\Omega^j$ is isomorphic to ${\rm SO}(d,1)\times_{{\rm
SO}(d),\xi_j}\wedge^j{\mathbb R}^d.$ By an inclusion:
\[
 {\rm SO}(d-1)={\rm SO}(2n)\to {\rm SO}(d),\quad A \mapsto
\left(\begin{array}{cc}
A&0\\
0&1
\end{array}
\right)
\]
the restriction of $\xi$ to ${\rm SO}(2n)$ is decomposed into a direct
sum of the standard representation $\sigma$ of ${\rm SO}(2n)$ on ${\mathbb
R}^{2n}$ and the trivial module ${\bf 1}$. Therefore we have
\[
 \xi_j|_{{\rm SO}(2n)}\simeq \sigma_j\oplus \sigma_{j-1},
\] 
where $\sigma_{j}$ is the $j$-th exterior product of $\sigma$. Let us observe that $\sigma_0$
and $\sigma_{2n}$ are trivial and $\sigma_j$ is isomorphic to
$\sigma_{2n-j}$. 
$\sigma_j\otimes{\mathbb C}$ is irreducible for $j\neq n$ whereas
$\sigma_n\otimes{\mathbb C}$ splits into a direct sum of two irreducible
representations, $\sigma_n^{+}$ and $\sigma_n^{-}$.
We prepare notation. Let
$\gamma\in \Gamma$ be a hyperbolic element. Then it is conjugate an
element of $MA$, $m_{\gamma} \exp[l(\gamma)H]$ $(m_{\gamma}\in M,
\,\exp[l(\gamma)H]\in A)$ where $l(\gamma)$ is the length of the
conjugacy class of $\gamma$. There is a $\gamma_0\in \Gamma$ which
determines a prime conjugacy class and that $\gamma=\gamma_0^{\mu(\gamma)}$,
where $\mu(\gamma)$ is a positive integer. For $0\leq j\leq 2n$ we put
\[
 \alpha_j(\gamma)=\frac{{\rm Tr}\rho(\gamma)\cdot {\rm
 Tr}\sigma_j(m_{\gamma})\cdot l(\gamma_0)}{\Delta(\gamma)},\quad c_j=|n-j|
\]
and
\[
 S_j(z)=\exp[-\sum_{\gamma\in \Gamma_{hyp}}\frac{\alpha_j(\gamma)}{l(\gamma)}e^{-zl(\gamma)}],
\]
where
\[
 \Delta(\gamma)=\det(I_{d-1}-e^{-l(\gamma)}m_{\gamma}).
\]
Let $s_j$ be the logarithmic derivative of $S_j$:
\[
 s_j(z)=\frac{d}{dz}\log S_j(z)=\sum_{\gamma\in \Gamma_{hyp}}\alpha_j(\gamma)e^{-zl(\gamma)}.
\]
If ${\rm Re}z$ is sufficiently large $S_j(z)$ absolutely converges. 

\begin{fact}(\cite{Fried}, (RS))
\[
 R_X(z,\rho)=\prod_{j=0}^{2n}S_j(z+j)^{(-1)^{j+1}}.
\]
\end{fact}
An isomorphism $\sigma_j\simeq \sigma_{2n-j}$ induces
\[
 \alpha_j(\gamma)=\alpha_{2n-j}(\gamma)\quad \mbox{and}\quad S_j(z)=S_{2n-j}(z),
\]
and therefore 

\begin{lm}
\begin{eqnarray*}
 r_X(z,\rho)&=&\sum_{j=0}^{2n}(-1)^{j+1}s_j(z+j)\\
&=& \sum_{j=0}^{n-1}(-1)^{j+1}\{s_j(z+j)+s_j(z+2n-j)\}+(-1)^{n+1}s_n(z+n).
\end{eqnarray*}
\end{lm}
\begin{lm}
Let $f$ be a meromorphic function defined on a neighborhood of the
 origin and
\[
 f(z)=a_0z^h(1+a_1z+\cdots), \quad a_0\neq 0
\]
its Laurent expansion. Then
\[
 a_1=\frac{1}{2}\lim_{z\to 0}\{\frac{f^{\prime}}{f}(z)+\frac{f^{\prime}}{f}(-z)\}.
\]
\end{lm}
{\bf Proof.} Set $f(z)=a_0z^hg(z)$, where $g(z)=1+a_1z+\cdots.$  Then
\[
 \frac{f^{\prime}}{f}(z)=\frac{h}{z}+\frac{g^{\prime}}{g}(z).
\]
and therefore 
\[
 \lim_{z\to 0}\{\frac{f^{\prime}}{f}(z)+\frac{f^{\prime}}{f}(-z)\}=\frac{2g^{\prime}(0)}{g(0)}=2a_1.
\]
\begin{flushright}
$\Box$
\end{flushright}
Let us regard a meromorphic continuation of $R_X(z,\rho)$ for a moment 
and
\[
 R_X(z,\rho)=c_0z^h(1+c_1z+\cdots),\quad c_0\neq 0
\]
the Taylor expansion at the origin. By {\bf Lemma 2.2} we obtain 
\[
 c_1=\frac{1}{2}\lim_{z\to 0}\{r_X(z,\rho)+r_X(-z,\rho)\}.
\]
Using Selberg trace formula we will compute RHS. Let $\Delta_X^j$ be
Hodge Laplacian acting on the space of smooth
sections of $\Omega^j_{X}(\rho)$ and its selfadjoint extension to
$L^2(X, \Omega^j_X(\rho))$ will be denoted by the same character. Since
$\rho$ is cuspidal $\Delta_X^j$ has only discrete spectrum which do not
accumulate and Selberg trace fomula for the heat kernel becomes
\[
 {\rm Tr}[e^{-t\Delta_X^j}]=H_j(t)+I_j(t)+U_j(t), \quad t>0.
\]
Here $H_j(t)$, $I_j(t)$ and $U_j(t)$ are the hyperbolic, the identical
and the unipotent orbital integral,
respectively(\cite{Sarnak-Wakayama}). In this section we will compute the derivative of Laplace transform of an each of
them:
\[
 L(f)(z)=2z\int^{\infty}_0e^{-tz^2}f(t)dt, \quad f=H_j,\,I_j,\,U_j.
\]

\subsection{The hyperbolic orbital integral}

For $0\leq j\leq 2n$, let us put
\[
 h_j(t)=\frac{1}{\sqrt{4\pi t}}\sum_{\gamma\in \Gamma_{hyp}}\alpha_j(\gamma)\exp\{-(\frac{l(\gamma)^2}{4t}+tc_j^2+nl(\gamma))\}.
\]
Then a hyperbolic orbital integral is given by (\cite{Fried}, {\bf
Theorem 2})
\[
 H_j(t)=h_j(t)+h_{j-1}(t),
\]
where $h_{-1}(t)$ is understood to be $0$. Notice that
\[
 h_j(t)=h_{2n-j}(t), \quad 0\leq j
\leq n
\]
Although the following lemma
seems to be well known, we will give a proof for a completeness.
\begin{lm}
Let $l$ and $z$ be positive numbers. Then
\[
 \int^{\infty}_0 \frac{1}{\sqrt{\pi t}}e^{-z^2t-\frac{l^2}{4t}}dt=\frac{e^{-lz}}{z}.
\]
\end{lm}
{\bf Proof.} Let us remember the well known formula:
\[
 \int^{\infty}_0 e^{-t^2-\frac{x^2}{t^2}}dt=\frac{\sqrt{\pi}}{2}e^{-2x},\quad x>0.
\]
If we differentiate it with respec to $x$, we obtain
\[
 x\int^{\infty}_0\frac{1}{t^2} e^{-t^2-\frac{x^2}{t^2}}dt=\frac{\sqrt{\pi}}{2}e^{-2x}.
\]
A change of variables, $t=z\sqrt{y},\, x=\frac{lz}{2}$ will yield
\[
 \frac{1}{\sqrt{4\pi}}\int^{\infty}_0 y^{-\frac{3}{2}}e^{-z^2y-\frac{l^2}{4y}}dy=\frac{e^{-lz}}{l}.
\]
Take a derivative of this equation with respect
to $z$, the desired formula will be proved. 
\begin{flushright}
$\Box$
\end{flushright}
Therefore if $z$ is a sufficiently large positive number,
\begin{eqnarray*}
L(e^{tc_j^2}h_j)(z) &=& 2z\int^{\infty}_0\sum_{\gamma\in
 \Gamma_{hyp}}\alpha_j(\gamma)e^{-nl(\gamma)}\frac{1}{\sqrt{4\pi
 t}}e^{-z^2t-\frac{l(\gamma)^2}{4t}}dt\\
&=&\sum_{\gamma\in
 \Gamma_{hyp}}\alpha_j(\gamma)e^{-nl(\gamma)}z\int^{\infty}_0\frac{1}{\sqrt{\pi
 t}}e^{-z^2t-\frac{l(\gamma)^2}{4t}}dt\\
&=& \sum_{\gamma\in
 \Gamma_{hyp}}\alpha_j(\gamma)e^{-(z+n)l(\gamma)}\\
&=& s_j(z+n).
\end{eqnarray*}
and we have proved the following proposition.
\begin{prop}
For a sufficiently large positive number $z$,
\begin{equation}
 L(e^{tc_j^2}h_j)(z)=s_j(z+n).
\end{equation}
\end{prop}
Park has obtained the following proposition even though $\rho$ is not
cuspidal (\cite{Park2007}). For the sake
of a convenience, we will give a proof in \S2.4 under our assumption.
\begin{prop}
$s_j(z)$ is continued to the entire plane as a meromorphic function
 whose singularities are at most only simple
 poles with integral residues. 
\end{prop}
This implies a meromorphic continuation of $S_j$ to the whole
plane. Therefore by {\bf Fact 2.1} $R_X(z,\rho)$ is also meromorphically
continued.
\subsection{The identical orbital integral}
We put 
\[
 i_j(t)=i_{2n-j}(t)=\frac{r}{4\pi}vol(\Gamma\backslash
 G)\int^{\infty}_{-\infty}e^{-t(\lambda^2+c_j^2)}P_j(\lambda)d\lambda,
 \quad 0\leq j\leq n-1,
\]
and
\[
 i_n(t)=\frac{r}{2\pi}vol(\Gamma\backslash
 G)\int^{\infty}_{-\infty}e^{-t\lambda^2}P_n(\lambda)d\lambda.
\]
Here $vol(\Gamma\backslash G)$ is the volume of $\Gamma\backslash G$ and $P_j$ is the Plancherel measure for $\sigma_j$ (\cite{Miatello}):
\[
 P_{j}(\lambda)=\frac{4^{1-n}}{(2n-1)!!^2\pi}
\left(\begin{array}{c}2n\\j\end{array}\right)q_j(\lambda),
\] 
where
\[
 q_j(\lambda)=\prod_{k=1}^j\{\lambda^2+(n-k+1)^2\}\prod_{k=j+1}^{n}\{\lambda^2+(n-k)^2\}
\]
Then the identical orbital integral is given by 
\[
 I_j(t)=i_j(t)+i_{j-1}(t).
\]
Since $\Gamma$ is torsin free its intersection with $K$ is the only
identity element. Remember that we have normalized the Haar measure so
that $vol(K)$ to be one and thus
\[
 vol(\Gamma\backslash G)=vol(\Gamma\backslash{\mathbb H}^d)\cdot vol(K)=vol(X).
\]
Hence we have obtained 
\[
 i_j(t)=i_{2n-j}(t)=\frac{r}{4\pi}vol(X)\int^{\infty}_{-\infty}e^{-t(\lambda^2+c_j^2)}P_j(\lambda)d\lambda,
 \quad 0\leq j \leq n-1,
\]
and
\[
 i_n(t)=\frac{r}{2\pi}vol(X)\int^{\infty}_{-\infty}e^{-t\lambda^2}P_n(\lambda)d\lambda.
\]
For example if $d=3$, i.e. $n=1$, using
\begin{equation}
 \int^{\infty}_{-\infty}e^{-t\lambda^2}d\lambda=\sqrt{\pi}t^{-\frac{1}{2}},
\end{equation}
one can see
\begin{equation}
 i_0(t)=i_{2}(t)=\frac{r}{4\pi^2}vol(X)\int^{\infty}_{-\infty}e^{-t(\lambda^2+1)}\lambda^2d\lambda=\frac{r\cdot
 vol(X)}{8\pi\sqrt{\pi}}e^{-t}t^{-\frac{3}{2}},
\end{equation}
\begin{equation}
 i_1(t)=\frac{r}{\pi^2}vol(X)\int^{\infty}_{-\infty}e^{-t\lambda^2}(\lambda^2+1)d\lambda=\frac{r\cdot
 vol(X)}{2\pi\sqrt{\pi}}(2t^{-\frac{1}{2}}+t^{-\frac{3}{2}}).
\end{equation}
In order to compute $L(e^{tc_j^2}i_j)(z)$ let us expand $q_j(\lambda)$ as
\[
 q_j(\lambda)=\sum_{k=0}^n\gamma_{j,k}\lambda^{2k},
\]
where $\gamma_{j,k}$ is an integer and $\gamma_{j,n}=1$.
\begin{lm} Let $z$ be a positive number. 
Then for a nonnegative integer $k$,
\[
 L(\int^{\infty}_{-\infty}e^{-t\lambda^2}\lambda^{2k}d\lambda)(z)=(-1)^k2\pi z^{2k}.
\]
In particular
 $L(\int^{\infty}_{-\infty}e^{-t\lambda^2}\lambda^{2k}d\lambda)(z)$ is
 entirely continued to the whole plane.
\end{lm}
{\bf Proof.} Let us take a $k$-times derivative of (4) with respect to $t$. Then we obtain
\[
 \int^{\infty}_{-\infty}e^{-t\lambda^2}\lambda^{2k}d\lambda=2^{-k}(2k-1)!!\sqrt{\pi}t^{-\frac{1}{2}-k},
\]
and
\[
 L(\int^{\infty}_{-\infty}e^{-t\lambda^2}\lambda^{2k}d\lambda)(z)=2^{-k}(2k-1)!!\sqrt{\pi}L(t^{-\frac{1}{2}-k})(z).
\]
Now since
\begin{eqnarray*}
L(t^{-\frac{1}{2}-k})(z) &=&
 2z\int^{\infty}_{0}e^{-tz^2}t^{-\frac{1}{2}-k}dt\\
&=& 2\Gamma(\frac{1}{2}-k)z^{2k}\\
&=& \frac{(-1)^k2^{k+1}\sqrt{\pi}}{(2k-1)!!}z^{2k},
\end{eqnarray*}
the desired equation has been proved.
\begin{flushright}
$\Box$
\end{flushright}

\begin{prop} We have 
\[
 L(e^{tc_j^2}i_j)(z)=L(e^{tc_{2n-j}^2}i_{2n-j})(z)=\frac{4^{1-n}r}{2(2n-1)!!^2\pi}
\left(\begin{array}{c}2n\\j\end{array}\right) vol(X)\sum_{k=0}^n(-1)^k\gamma_{j,k}z^{2k}.
\]
for $0 \leq j \leq n-1$ and
\[
 L(i_n)(z)=\frac{4^{1-n}r}{(2n-1)!!^2\pi}
\left(\begin{array}{c}2n\\n\end{array}\right) vol(X)\sum_{k=0}^n(-1)^k\gamma_{n,k}z^{2k}.
\]
\end{prop}
\begin{cor}
Suppose $n=1$. Then
\[
 L(e^ti_0)(z)= L(e^ti_2)(z)=-\frac{r}{2\pi}vol(X)z^2
\]
and
\[
 L(i_1)(z)=\frac{2r}{\pi}vol(X)(1-z^2).
\]
\end{cor}
\subsection{The unipotent orbital integral}
Let $\zeta_{\nu}(s,\chi_{\nu,i})$ be the Epstein L-function:
\[
 \zeta_{\nu}(s,\chi_{\nu,i})=\sum_{0\neq \eta\in \Gamma_{\nu}}\chi_{\nu,i}(\eta)|X_{\eta}|^{-2n(s+1)},
\]
where $X_{\eta}$ is an element of the Lie algebra of $N_{\nu}$ such that
$\exp(X_{\eta})=\eta$. The norm is taken with respect to the
normalized Cartan-Killng form. It absolutely converges if ${\rm Re}\,s$ is
sufficiently large and is meromorphically continued to the whole
plane. Since $\chi_{\nu,i}$ is nontrivial it is regular at the
origin and we put
\[
 \tau_{\nu}=\sum_{i=1}^r\zeta_{\nu}(0,\chi_{\nu,i}).
\] 
Let $A(\frak{n})$ be the volume of the unit spehre in ${\frak n}$. By
\cite{Osborne-Warner} we find the unipotent orbital integral is given by 
\[
 U_j(t)=u_j(t)+u_{j-1}(t),
\]
where
\begin{eqnarray}
 u_j(t)=u_{2n-j}(t)&=&\frac{1}{2\pi A({\frak
 n})}\sum_{\nu=1}^hvol(\Gamma_{\nu}\backslash
 N_{\nu})\tau_{\nu}\int^{\infty}_{-\infty}e^{-t(\lambda^2+c_j^2)}d\lambda\\
&=& \frac{1}{2\sqrt{\pi} A({\frak
 n})}\sum_{\nu=1}^hvol(\Gamma_{\nu}\backslash N_{\nu})\tau_{\nu}e^{-tc_j^2}t^{-\frac{1}{2}}.
\end{eqnarray}
for $0\leq j\leq n-1$ and
\begin{eqnarray}
 u_n(t)&=&\frac{1}{\pi A({\frak
 n})}\sum_{\nu=1}^hvol(\Gamma_{\nu}\backslash
 N_{\nu})\tau_{\nu}\int^{\infty}_{-\infty}e^{-t\lambda^2}d\lambda\\
&=& \frac{1}{\sqrt{\pi} A({\frak
 n})}\sum_{\nu=1}^hvol(\Gamma_{\nu}\backslash N_{\nu})\tau_{\nu}t^{-\frac{1}{2}}.
\end{eqnarray}
Here we have used (4).
\begin{lm}Let $z$ be a positive number. Then
\[
 L(\int^{\infty}_{-\infty}e^{-t\lambda^2}d\lambda)(z)=2\pi,
\]
and $L(\int^{\infty}_{-\infty}e^{-t\lambda^2}d\lambda)(z)$ is entirely
 continued to the whole plane as $2\pi$.
\end{lm}
{\bf Proof.} We compute,
\begin{eqnarray*}
L(\int^{\infty}_{-\infty}e^{-t\lambda^2}d\lambda)(z) &=& 2z
 \int^{\infty}_0dt
 e^{-tz^2}\int^{\infty}_{-\infty}e^{-t\lambda^2}d\lambda\\
&=& 2z \int^{\infty}_{-\infty}d\lambda  \int^{\infty}_0
 e^{-t(\lambda^2+z^2)}dt\\
&=& 2z\int^{\infty}_{-\infty}\frac{d\lambda}{\lambda^2+z^2}.
\end{eqnarray*}
The desired formula will be obtained by the contour integration.
\begin{flushright}
$\Box$
\end{flushright}
Thus putting 
\[
 \delta(X,\rho)=\frac{1}{A({\frak n})}\sum_{\nu=1}^hvol(\Gamma_{\nu}\backslash N_{\nu})\tau_{\nu},
\]
we have proved the following proposition.
\begin{prop} For $0\leq j\leq n-1$, both $L(e^{tc_j^2}u_j)(z)$ and
 $L(e^{tc_{2n-j}^2}u_{2n-j})(z)$ are analytically continued to the entire plane as a
 constant $\delta(X,\rho)$, whereas $L(u_n)(z)$ is continued as
$2\delta(X,\rho).$
\end{prop}

\subsection{An application of Selberg trace formula}
Let $0 \leq j \leq n$. Then by definition we know
\[
 h_j(t)=\sum_{k=0}^j(-1)^{j-k}H_k(t), \quad i_j(t)=\sum_{k=0}^j(-1)^{j-k}I_k(t)
\]
and
\[
 u_j(t)=\sum_{k=0}^j(-1)^{j-k}U_k(t).
\]
If we put
\[
 \delta_j(t)=\sum_{k=0}^j(-1)^{j-k}{\rm Tr}[e^{-t\Delta^k_X}],
\]
Selberg trace formula implies
\begin{equation}
\delta_j(t)=h_j(t)+i_j(t)+u_j(t).
\end{equation}
The following lemma will directly follow from the definition of $L$.
\begin{lm}
\[
 L(e^{tc_j^2}\delta_j)(-z)=-L(e^{tc_j^2}\delta_j)(z).
\]
\end{lm}

Now we will prove {\bf Proposition 2.2}. Since
$\rho$ is cuspidal, $\Delta_X^k$ has only discrete spectrum $\{\sigma_k(l)\}_l$ which do not
accumulate and are nonnegative. Let us fix
$z\in{\mathbb C}$ so that
\[
 {\rm Re}z^2 > n^2. 
\]
Then
\begin{eqnarray*}
L(e^{tc_j^2}\delta_j)(z) &=& \sum_{k=0}^j(-1)^{j-k}\sum_{l}2z\int^{\infty}_0e^{-t(z^2-c_j^2+\sigma_k(l))}dt\\
&=& \sum_{k=0}^j(-1)^{j-k}\sum_{l}\frac{2z}{z^2-c_j^2+\sigma_k(l)},
\end{eqnarray*}
and 
\begin{equation}
 L(e^{tc_j^2}\delta_j)(z-c_j)= \sum_{k=0}^j(-1)^{j-k}\sum_{l}\{\frac{1}{z-c_j+\sqrt{c_j^2-\sigma_k(l)}}+\frac{1}{z-c_j-\sqrt{c_j^2-\sigma_k(l)}}\}.
\end{equation}
Thus $L(e^{tc_j^2}\delta_j)(z)$ is meromorphically continued to the
whole plane and has
only simple poles with integral residues. Using {\bf Proposition 2.1},
{\bf Proposition 2.3} and {\bf Proposition 2.4} Selberg
trace formula imply
\begin{eqnarray}
 s_j(z+n)&=& L(e^{tc_j^2}\delta_j)(z)\nonumber\\
&-&\frac{4^{1-n}r(1+\delta_{j,n})}{2(2n-1)!!^2\pi}\left(\begin{array}{c}2n\\j\end{array}\right)vol(X)\sum_{k=0}^n(-1)^k\gamma_{j,k}z^{2k}\nonumber\\
&-& (1+\delta_{j,n})\delta(X,\rho),
\end{eqnarray}
where $\delta_{j,n}$ is the Kronecker's delta. This proves {\bf
Proposition 2.2}.
\begin{flushright}
$\Box$
\end{flushright}
Using {\bf Lemma 2.1} and (13), the above computation shows
\begin{eqnarray*}
r_X(z,\rho) &=&
 \sum_{j=0}^{n-1}(-1)^{j+1}\sum_{k=0}^j(-1)^{j-k}\sum_{l}\{\frac{1}{z-c_j+\sqrt{c_j^2-\sigma_k(l)}}+\frac{1}{z-c_j-\sqrt{c_j^2-\sigma_k(l)}}\}\\
&+&
 \sum_{j=0}^{n-1}(-1)^{j+1}\sum_{k=0}^j(-1)^{j-k}\sum_{l}\{\frac{1}{z+c_j+\sqrt{c_j^2-\sigma_k(l)}}+\frac{1}{z+c_j-\sqrt{c_j^2-\sigma_k(l)}}\}\\
&+&
 (-1)^{n+1}\sum_{k=0}^n(-1)^{n-k}\sum_{l}\{\frac{1}{z+\sigma_k(l)\sqrt{-1}}+\{\frac{1}{z-\sigma_k(l) \sqrt{-1}}\}\\
&+& E(z),
\end{eqnarray*}
where $E(z)$ is an entire function. Thus remembering $h^0(X,\rho)=0$, we
have 
\[
 {\rm Res}_{z=0}r_X(z,\rho)=2\sum_{l=0}^{n-1}(-1)^l(n-l)h^{l+1}(X,\rho).
\]
It is obvious that the zeros and poles of $R_X(z,\rho)$ are located on
\[
 \Xi=\{z\in{\mathbb C}\,|\, {\rm Re}z=-n,-(n-1),\cdots, n-1,n\},
\]
except for finitely many of them. For example if $d=3$ (i.e. $n=1$), we
conclude $R_X(z,\rho)$ has a zero at the origin of order
$2h^1(X,\rho)$. If the minimum of spectrum of $\Delta_X^0$ is greater
than or equal to 1, all zeros and poles are located on $\Xi$. The
following theorem is a consequence of {\bf Lemma 2.1}, {\bf Lemma 2.6}
and (13).
\begin{thm}
\[
 r_X(z,\rho)+r_X(-z,\rho)=\frac{4^{1-n}r\cdot vol(X)}{(2n-1)!!^2\pi}\chi(z)+4\sum_{j=0}^n(-1)^j\delta(X,\rho),
\]
where
\[
 \chi(z)= \sum_{j=0}^{n}(-1)^{j}\left(\begin{array}{c}2n\\j\end{array}\right)\sum_{k=0}^n(-1)^k\gamma_{j,k}\{(z+j-n)^{2k}+(z-j+n)^{2k}\}.
\]
\end{thm}
\begin{cor} Suppose $d=3$. Then
\[
 r_X(z,\rho)+r_X(-z,\rho)=\frac{2r}{\pi} vol(X)(z^2-3).
\] 
\end{cor}

Noting the order of $R_X(z,\rho)$ is even, an easy computation will show
\[
 \lim_{z\to 0}R_X(z,\rho)R_X(-z,\rho)^{-1}=1.
\]
By {\bf Theorem 2.1}, 
\begin{eqnarray*}
\frac{d}{dz}\log(R_X(z,\rho)R_X(-z,\rho)^{-1}) &=&
 r_X(z,\rho)+r_X(-z,\rho)\\
&=& \frac{4^{1-n}r\cdot vol(X)}{(2n-1)!!^2\pi}\chi(z)+4\sum_{j=0}^n(-1)^n\delta(X,\rho),
\end{eqnarray*}
and therefore $R_X(z,\rho)$ satisfies a functional equation:
\[
 R_X(z,\rho)\cdot R_X(-z,\rho)^{-1}=\exp[\frac{4^{1-n}r\cdot vol(X)}{(2n-1)!!^2\pi}X(z)+4\sum_{j=0}^n(-1)^n\delta(X,\rho)z],
\]
where $X(z)$ is the primitive function of $\chi(z)$ so that
$X(0)=0$. For example if $d=3$,
\[
 R_X(z,\rho)\cdot R_X(-z,\rho)^{-1}=\exp[\frac{2r}{\pi} vol(X)(\frac{z^3}{3}-3z)].
\]
Now the second coefficient of the Taylor expansion is obtained by {\bf
Lemma 2.2} and {\bf Theorem 2.1}.
\begin{thm}
Let 
\[
 R_X(z,\rho)=c_0z^h(1+c_1z+\cdots), \quad c_0\neq 0,
\]
be the Laurent expansion. Then $c_1-2\sum_{j=0}^n(-1)^n\delta(X,\rho)$ is a rational
 multiple of $vol(X)/\pi$. 
\end{thm}
\begin{cor}Suppose $d=3$. Then
\[
 c_1=-\frac{3r}{\pi} vol(X).
\]
\end{cor}
 
\section{The leading coefficient}
Throughout this section we assume that $X$ is a hyperbolic threefold with
finite volume. We will compute the leading coefficient of the Taylor
expansion of $R_X(z,\rho)$ at the origin. In \S2.4 we have seen ${\rm ord}_{z=0}R_X(z,\rho)=2h^1(X,\rho)$. The
following fact is a special case of \cite{Park2007}. 
\begin{fact}
\[
 \lim_{z\to 0}z^{-2h^{1}(X,\rho)}R_{X}(z,\rho)={\rm exp}(-\zeta_{X}^{\prime}(0,\rho)).
\]
\end{fact}
Here 
\[
  \zeta_{X}(z,\rho)=\sum_{p=0}^{3}(-1)^p p\cdot \zeta_{X}^{(p)}(z,\rho),
\]
where
\[
 \zeta_{X}^{(p)}(z,\rho)=\frac{1}{\Gamma(z)}\int^{\infty}_{0}\{{\rm Tr}[e^{-t\Delta^{p}_X}]-h^{p}(X,\rho)\}t^{z-1}dt.
\]
$\zeta_{X}^{(p)}(z,\rho)$ absolutely converges if ${\rm Re}z$ is
sufficiently large and is meromorphically continued to the whole
plane. In fact let us put
\[
 \theta_p(t)={\rm Tr}[e^{-t\Delta^{p}_X}]-h^{p}(X,\rho).
\]
Then the computation of orbital integrals in \S2 and Selberg trace
formula show that it has an asymptotic expansion on $(0,1]$ such that 
\begin{equation}
 \theta_p(t)\sim t^{-\frac{3}{2}}\sum_{l=0}^Nc_lt^l+O(t^{N-\frac{3}{2}}).
\end{equation}
Therefore if ${\rm Re}z > N-3/2$,
\[
  \int^{1}_0 \theta_p(t)t^{z-1}dt=\sum_{l=0}^N\frac{c_l}{z+l-\frac{3}{2}}+R_N(z),
\]
where $R_N(z)$ is a regular function on $\{z\in{\mathbb
C}\,|\, {\rm Re}z > \frac{3}{2}-N\}$ which is meromorphically continued
to the whole plane. Since $\theta_p(t)$ exponentially decays as $t\to
\infty$, $\int^{\infty}_1 \theta_p(t)t^{z-1}dt$ is an entire
function. Thus writing
\[
 \int^{\infty}_0 \theta_p(t)t^{z-1}dt=\int^{1}_0 \theta_p(t)t^{z-1}dt+\int^{\infty}_1 \theta_p(t)t^{z-1}dt,
\]
we know that $\zeta_{X}^{(p)}(z,\rho)$ is meromorphically continued to
the whole plane and that it vanishes at the origin. Since we assume that
$\rho$ is cuspidal it is possible to prove {\bf Fact 3.1} just following
the arguments of \cite{Fried}. For a convenience we will give a proof in
{\bf Appendix}. Thus the leading coefficient is ${\rm
exp}(-\zeta_{X}^{\prime}(0,\rho))$ but we want to express this by a more
geometric term.
\subsection{Boundary conditions}

We will use the Poincar\'{e} upper half space model:
\[
 {\mathbb H}^3=\{(x,y,r)\in{\mathbb
R}^3\,|\,r>0\}, \quad g=\frac{dx^2+dy^2+dr^2}{r^2}.
\]
For $a\in{\mathbb R}$ we put
\[
 {\mathbb H}^{3}_{a}=\cap_{\nu=1}^{h}g_{\nu}{\mathbb H}^{3}_{a,\infty},
 \quad {\mathbb H}^{3}_{a,\infty}=\{(x,y,r)\in{\mathbb H}^3\,|\, r \leq e^{a}\}.
\]
(Remember that $g_{\nu}\in {\rm PSL}_2({\mathbb C})$ is chosen to
satisfy
\[
 N_{\nu}=g_{\nu}Ng_{\nu}^{-1}, \quad N=\{\left(
\begin{array}{cc}
1&z\\
0&1\end{array}
\right)\,|\,z\in{\mathbb C}\}.)
\]
Let $X_a$ be the image of ${\mathbb H}^{3}_{a}$ by the natural
projection
\[
 {\mathbb H}^{3}\stackrel{\pi}\to X,
\]
and $Y_a$ the closure of $X\setminus X_a$. If $a$ is sufficiently large
$Y_a$ is a disjoint union of $Y_{a,\nu}$ ($1 \leq \nu \leq h$) which is
a warped product of a flat 2-torus $T_{\nu}=N_{\nu}\slash
\Gamma_{\nu}$($={\mathbb R}^2/\Gamma_{\nu}$) and an interval $[e^a,\,\infty)$ with the metric,
\[
 g=du^2+e^{-2u}(dx^2+dy^2).
\]
(Here we have made a change of variables: $r=e^u$.) According to the
decomposition
\[
 \rho_{\nu}=\oplus_{i=1}^{r}\chi_{\nu,i},
\] 
a smooth section $\varphi$ of $\Omega^{p}_{X}(\rho)$ over $Y_{a,\nu}$ is
written as 
\[
 \varphi=\sum_{i=1}^{r}\varphi_i,\quad \varphi_i=\sum_{|\alpha|=p}\varphi_{i,\alpha}dx^{\alpha}\in
 C^{\infty}(Y_{a,\nu}, \Omega_{X}^{p}(\chi_{\nu,i})),
\] 
where we have put
\[
 x_0=u,\quad x_1=x,\quad x_2=y.
\]
\begin{lm}
\[
 \int_{T_{\nu}}\varphi_{i,\alpha}dxdy=0.
\]
\end{lm}
{\bf Proof.} Let us choose $\gamma\in \Gamma_{\nu}$ so that
\[
 \chi_{\nu,i}(\gamma)\neq 1.
\]
By definition we have
\[
 \gamma^{*}\varphi_{i,\alpha}=\chi_{\nu,i}(\gamma)\varphi_{i,\alpha}.
\]
and the desired result will follow from
\[
 \int_{T_{\nu}}\varphi_{i,\alpha}dxdy=\int_{T_{\nu}}\gamma^{*}\varphi_{i,\alpha}dxdy=\chi_{\nu,i}(\gamma)\int_{T_{\nu}}\varphi_{i,\alpha}dxdy.
\]
\begin{flushright}
$\Box$
\end{flushright}

We will consider an eigenvalue problem of Hodge Laplacian $\Delta^p_X$ on
spaces $L^2(X_a,\Omega^p_{X}(\rho))$ or
$L^2(Y_{a,\nu},\Omega^p_{X}(\rho))$ under a certain boundary
condition. Hereafter for simplicity we omit the subscript $X$ of $\Delta^p_X$. 
The restriction $\Omega_{X}^p(\rho)$ to the boundary
$T_{\nu}$ of $Y_{a,\nu}$ is decomposed into
\[
 \Omega_{X}^p(\rho)|_{T_{\nu}}=\Omega^{p}_{T_{\nu}}(\rho)\oplus du\wedge \Omega^{p-1}_{T_{\nu}}(\rho),
\]
and according to this a section $\omega$ of
$\Omega^p_{X}(\rho)|_{T_{\nu}}$ is expressed by
\[
 \omega=\omega_{tan}+ \omega_{norm},
\]
where $\omega_{tan}$ (resp. $\omega_{norm}$) is a section of
$\Omega^{p}_{T_{\nu}}(\rho)$ (resp. $du\wedge\Omega^{p-1}_{T_{\nu}}(\rho)$)

\begin{df} Let $\omega$ be a smooth section of $\Omega^p_{X}(\rho)$ on
 $X_a$ or $Y_{a,\nu}$. We call it satisfies {\rm the absolute boundary condition} if
 both $\omega_{norm}$ and $(d\omega)_{norm}$ vanish on the every
 boundary. If the Hodge dual $*\omega$ satisfies the absolute boundary
 condition we will refer that it satisfies {\rm the relative boundary
 condition}. If $\omega$ satisfies both of the absolute
 and the relative boundary condition, we call it satisfies {\rm the Dirichlet
 boundary condition}.
\end{df}
It is easy to see that $\omega$ satisfies the relative boundary
condition if and only if both $\omega_{tan}$ and $(d\omega)_{tan}$ vanish on every
 $T_{\nu}$.
Thus $\omega$ satisfies the Dirichlet boundary condition if and only if
the restrictions of both $\omega$
and $d\omega$ become the $0$-section of $\Omega_{X}^p(\rho)|_{T_{\nu}}$
and $\Omega_{X}^{p+1}(\rho)|_{T_{\nu}}$ for every $1\leq\nu\leq h$. More
concretely  the latter
condition means that if we write
\[
 \omega=\sum_{|\alpha|=p}f_{\alpha}dx^{\alpha}, \quad d\omega=\sum_{|\beta|=p+1}g_{\beta}dx^{\beta},
\]
all $f_{\alpha}$ and $g_{\beta}$ vanish along $T_{\nu}$ for every $\nu$.
Notice that $*$ interchanges the absolute and relative boundary conditions and
preserves the Dirichlet one. Since $\rho$ is unitary the associated local system possesses a fiberwise
hermitian inner product ${\rm Tr}_{\rho}$. For $\omega,\,\eta\in
\Omega^p_{X}(\rho)$ we put
\[
 (\omega,\eta)=\frac{{\rm Tr}_{\rho}(\omega\wedge *\eta)}{dv_g},
\]
which defines a hermitian inner product on $\Omega_{X}^p(\rho)$. Here
$dv_g$ is the volume form of $g$. Let $M$ be $X_{a}$ or $Y_{a,\nu}$ and
$\nabla$ the covariant derivative. If both $\omega$ and $\eta$
satisfy one of the boundary conditions,
\[
 \int_{M}(\Delta^p\omega,\eta)dv_g=\int_M(\nabla\omega,\nabla\eta)dv_g=\int_{M}(\omega,\Delta^p\eta)dv_g,
\]
by Stokes theorem. Therefore $\Delta^p$ has a
selfadjoint extension $\Delta^p_{abs}$, $\Delta^p_{rel}$ or
$\Delta^p_{dir}$ according to a boundary condition. If $*$ is $abs$ (resp.
$rel$ or $dir$) {\it its dual} $\hat{*}$ is defined to be $rel$ (resp.
$abs$ or $dir$). Since the Hodge $*$-operator intertwines the action of $\Delta^p_{*}$ on
$L^2(M,\Omega^p_X({\rho}))$ and one of $\Delta^{3-p}_{\hat{*}}$
on $L^2(M,\Omega^{3-p}_X({\rho}))$, we
will only consider the case  of $p=0$ or $1$. \\

For a later purpose we will
introduce one more boundary condition. Let $\alpha$ be a real number greater than
one. For a sufficiently large $a$,
$Y_{a,\nu}\cap X_{\alpha a}$ is diffeomorphic to
 $T_{\nu}\times[e^a,\alpha e^a]$.  If
$\omega\in C^{\infty}(Y_{a,\nu}\cap X_{\alpha a}, \Omega^p_X({\rho}))$ satisfies
 the Dirichlet condition on $T_{\nu}\times\{e^a\}$ and a condition $*$
 ($*=abs$, $rel$ or $dir$) on
$T_{\nu}\times\{\alpha e^a\}$ we call it satisfies {\it
Dirichlet/$*$-condition}.
If $\omega\in C^{\infty}(Y_{a}\cap X_{\alpha a},
\Omega^p_X({\rho}))$ satisfies Dirichlet/$*$-condition on every
connected component ($*$ does not depend on a component) we will refer that it satisfies
Dirichlet/$*$-condition.

\subsection{Spectrum of Hodge Laplacian at cusps}
Since $\Delta^p$ commutes with the
action of $\Gamma$ it preserves the decomposition,
\[
 C^{\infty}(Y_{a,\nu},\Omega^{p}_X(\rho)) =\oplus_{i=1}^{r}C^{\infty}(Y_{a,\nu},\Omega^{p}_X(\chi_{\nu,i})).
\]
Thus the spectral problem of Hodge Laplacian on
$L^2(Y_{a,\nu},\Omega^{p}_X(\rho))$ is reduced to one on
$L^{2}(Y_{a,\nu},\Omega^{p}_X(\chi_{\nu,i}))$. We will give an explicit
formula of $\Delta^p$ on $Y_{a,\nu}$. A straightforward  computation will show the following lemma.
\begin{lm} Let $\Delta_{T}$ be the positive Laplacian on a flat torus,
\[
 \Delta_T=-(\partial_x^2+\partial_y^2).
\]
\begin{enumerate}
\item For $f\in C^{\infty}(Y_{a,\nu},\Omega^0_X(\chi_{\nu,i}))$,
\[
 \Delta^0f=e^{2u}\Delta_Tf-\partial_u^2f+2\partial_uf.
\]
\item For $\omega=fdx+gdy+hdu\in
      C^{\infty}(Y_{a,\nu},\Omega^1_X(\chi_{\nu,i}))$,
\begin{eqnarray*}
 \Delta^1\omega &=& (e^{2u}\Delta_Tf-\partial_u^2f+2\partial_xh)dx\\
&+& (e^{2u}\Delta_Tg-\partial_u^2g+2\partial_yh)dy\\
&+& (e^{2u}\Delta_Th-\partial_u^2h+2\partial_uh-2e^{2u}(\partial_xf+\partial_yg))du.
\end{eqnarray*}
\end{enumerate}

\end{lm}

\begin{fact}(\cite{Reed-Simon}{\bf Theorem XIII.1}, The min-max
 principle) Let $A$ be a selfadjoint operator with domain $D(A)$, which
 is bounded below. Define
\[
 \mu_n(A)=\sup_{\varphi_1,\cdots,\varphi_{n-1}}U_A(\varphi_1,\cdots,\varphi_{n-1}),
\] 
where
\[
 U_A(\varphi_1,\cdots,\varphi_{n-1})=\inf_{\psi\in
 D(A),||\psi||=1, \psi\in <\varphi_1,\cdots,\varphi_{n-1}>^{\perp}}(\psi,A\psi),
\]
and $<\varphi_1,\cdots,\varphi_{n-1}>^{\perp}$ is the orthogonal
 complement of a vector space $<\varphi_1,\cdots,\varphi_{n-1}>$ spanned
 by $\{\varphi_1,\cdots,\varphi_{n-1}\}$. Then either of the followings holds:
\begin{enumerate}
\item there are $n$ eingenvalues below the bottom of the essential
      spectrum and $\mu_n(A)$ is the $n$-th eigenvalue counting with
      multiplicity,
\item $\mu_n(A)$ is the bottom of the essential spectrum.
\end{enumerate}
\end{fact}

Later on we will need a variant of this. 
\begin{lm}
Let $A$ be a selfadjoint operator bounded below such that
 $(A-\lambda)^{-1}$ is compact for a certain $\lambda\in \rho(A)$, where
 $\rho(A)$ is the resolvent set. Then the $n$-th eigenvalue $\mu_n(A)$ is obtained by
\[
 \mu_n(A)=\inf_{\frak{M}\in {\rm Gr}_nD(A)}\sup_{0\neq v\in\frak{M}}\frac{(Av,v)}{||v||^2}.
\]
Here ${\rm Gr}_nD(A)$ is the set of $n$-dimensional subspaces of $D(A)$.
\end{lm}
{\bf Proof.}
Let $\mu^{\prime}_n(A)$ be rhe RHS of the above equation. By the
assumption there is a complete orthonormal basis $\{\varphi_n\}_n$ in
$D(A)$ such that $A\varphi_n=\mu_n(A)\varphi_n$ with $\mu_1(A) \leq
\mu_2(A)\leq \cdots$ and $\mu_n(A)\to \infty.$ Let $\frak{N}$ be an
$n$-dimensional space spanned by $\{\varphi_1,\cdots, \varphi_n\}$. Thus
\[
 \mu_n(A)=\sup_{0\neq v\in\frak{N}}\frac{(Av,v)}{||v||^2},
\]
and $\mu^{\prime}_n(A)\leq \mu_n(A)$ by definition. Suppose $\mu^{\prime}_n(A)$ is
strictly less than $\mu_n(A)$. Then there is an $n$-dimensional subspace
$\frak{M}$ of $D(A)$ so that
\[
 \mu^{\prime}_n(A) \leq \sup_{0\neq v\in\frak{M}}\frac{(Av,v)}{||v||^2}
 < \mu_n(A).
\] 
But by the equation (2a) in pp.77 of \cite{Reed-Simon}, the dimension of $\frak{M}$ should be less than $n$, which is a contradiction.
\begin{flushright}
$\Box$
\end{flushright}

Let $A$ and $B$ are selfadjoint operators bounded below which act on a Hilbert space $H$. Suppose that they have the same domain $D$ and
that $A\geq B$, i.e. $(Av,v)\geq (Bv,v)$ for any $v\in D$. Then {\bf
Fact 3.1} implies

\begin{lm}
\[
 \mu_n(A)\geq \mu_n(B)
\]
\end{lm}

Let $a$ and $a^{\prime}$ be positive numbers so that $a^{\prime}\geq
a$. Extending as $0$-map on the outside $L^{2}(X_a,\Omega^p_X(\rho))$ is
embedded into $L^{2}(X_{a^{\prime}},\Omega^p_X(\rho))$. Thus
$D(\Delta^p_{dir}|_{X_a})$ is a subspace of
$D(\Delta^p_{dir}|_{X_{a^{\prime}}})$. In particular ${\rm
Gr}_n(D(\Delta^p_{dir}|_{X_a}))$ is a subset of ${\rm
Gr}_n(D(\Delta^p_{dir}|_{X_{a^{\prime}}}))$. Since
$\Delta_{dir}^p|_{X_{a^{\prime}}}$ satisfies the assumption of {\bf
Lemma 3.3},
\[
 \mu_{n}(\Delta_{dir}^p|_{X_{a^{\prime}}}) \leq
 \mu_{n}(\Delta_{dir}^p|_{X_{a}}).
\]
Since $\rho$ is cuspidal $\Delta^p$ also satisfies the assumption of {\bf
Lemma 3.3}. The same argument will yield the following lemma.

\begin{lm}
\begin{enumerate}
\item Let $a$ and $a^{\prime}$ be positive numbers so that $a^{\prime}\geq a$. Then,
\[
 \mu_{n}(\Delta_{dir}^p|_{X_{a^{\prime}}}) \leq \mu_{n}(\Delta_{dir}^p|_{X_{a}}).
\]
\item For a positive $a$,
\[
 \mu_{n}(\Delta^p) \leq \mu_{n}(\Delta_{dir}^p|_{X_{a}}).
\]
and
\[
 \mu_{n}(\Delta_{*}^p|_{X_{a}}) \leq \mu_{n}(\Delta_{dir}^p|_{X_{a}}),
\]
where $*$ is $abs$ or $rel$.
\end{enumerate}
\end{lm}
\begin{remark}The above lemma also follows from the Rayleigh-Ritz
 technique.(\cite{Reed-Simon}{\bf Theorem XIII.3})
\end{remark}

Let $\Gamma_{\nu}^{*}$ be the dual lattice of $\Gamma_{\nu}$. We
will define its {\it norm} to be
\[
 ||\Gamma_{\nu}^{*}||={\rm Min}\{|\gamma|\,|\,0\neq \gamma\in\Gamma_{\nu}^{*}\}.
\]
Here the modulus $|\cdot|$ is taken with respect to the standard
Euclidean metric $dx^2+dy^2$ on ${\mathbb R}^2$.

\begin{prop} 
\[
 \mu_{1}(\Delta_{dir}^0|_{Y_{a,\nu}}) \geq e^{2a}||\Gamma_{\nu}^{*}||^2.
\]
\end{prop}
{\bf Proof.} Let us consider a nonnegative selfadjoint operator
\[
 P_a=e^{2a}\Delta_T-\partial_u^2+2\partial_u
\]
on $L^2(Y_{a,\nu},\Omega^0(\chi_{\nu,i}))$ under Dirichlet condition at
the boundary. Since 
\[
 \Delta^0-P_a=(e^{2u}-e^{2a})\Delta_T
\]
is a nonnegative operator {\bf Lemma 3.4} implies
\[
 \mu_1(\Delta^0_{dir}|_{Y_{a,\nu}})\geq \mu_1(P_a).
\]
For $f\in
C^{\infty}_c(Y_{a,\nu},\Omega^0(\chi_{\nu,i}))$,
\begin{eqnarray*}
\int_{Y_{a,\nu}}(P_af,f)dv_g &=&
 e^{2a}\int_{Y_{a,\nu}}\Delta_Tf\cdot\bar{f}e^{-2u}dxdydu +
 \int_{Y_{a,\nu}}|\partial_uf|^2e^{-2u}dxdydu\\
& \geq & e^{2a}\int_{Y_{a,\nu}}\Delta_Tf\cdot\bar{f}e^{-2u}dxdydu\\
&=& e^{2a}\int^{\infty}_{a}due^{-2u}\int_{T_{\nu}}\Delta_Tf\cdot\bar{f}dxdy.
\end{eqnarray*}
Let
\[
 f=\sum_{\gamma\in{\Gamma_{\nu}^{*}}}\{f_{\gamma}(u){\bf
 e}_{\gamma}(z)+f_{\gamma}^{*}(u){\bf e}_{\gamma}(\bar{z})\}, \quad {\bf
 e}_{\gamma}(z)=\exp (2\pi i\gamma z)
\]
be a Fourier expansion with respect to $T_{\nu}$-direcrion. Here notice
that by {\bf Lemma 3.1} $\gamma$ runs through nonzero
elements of $\Gamma_{\nu}^{*}$. Then
\begin{eqnarray*}
\int_{T_{\nu}}\Delta_Tf\cdot\bar{f} &=& {\rm vol}(T_{\nu})\sum_{0\neq
 \gamma \in
 \Gamma_{\nu}^{*}}|\gamma|^2\{|f_{\gamma}(u)|^2+|f_{\gamma}^{*}(u)|^2\}\\
&\geq & ||\Gamma_{\nu}^{*}||^2 {\rm vol}(T_{\nu}) \sum_{0\neq
 \gamma \in
 \Gamma_{\nu}^{*}}\{|f_{\gamma}(u)|^2+|f_{\gamma}^{*}(u)|^2\}\\
&=& ||\Gamma_{\nu}^{*}||^2\int_{T_{\nu}}|f|^2dxdy,
\end{eqnarray*}
and therefore we have obtained
\[
 \int_{Y_{a,\nu}}(P_af,f)dv_g \geq e^{2a}||\Gamma_{\nu}^{*}||^2\int_{Y_{a,\nu}}(f,f)dv_g.
\]
Now {\bf Fact 3.3} implies $\mu_1(P_a)\geq
e^{2a}||\Gamma_{\nu}^{*}||^2$.
\begin{flushright}
$\Box$
\end{flushright}
The same argument will prove
\begin{prop} For $\alpha>1$ and $*=abs$ or $rel$, 
\[
 \mu_{1}(\Delta_{dir\slash *}^0|_{X_{\alpha a}\cap Y_{a,\nu}}) \geq e^{2a}||\Gamma_{\nu}^{*}||^2.
\]
\end{prop}
Next we will estimate $\mu_{1}(\Delta_{dir}^1|_{Y_{a,\nu}})$ from
below.
Before doing this we will give some remarks. Let us fix a positive number
$\alpha$ less than $a$ and we make a change of variables,
\[
 u=v+\alpha.
\]
Then $Y_{a,\nu}$ is isometric to a warped product,
\[
 [a^{\prime}, \infty)\times T^{\prime}_{\nu},\quad a^{\prime}=a-\alpha,
\]
with metric 
\[
 dg=dv^2+e^{-2v}(dx^2+dy^2).
\]
Here the boundary $T^{\prime}_{\nu}$ is a quotient of ${\mathbb R}^2$ with
the standard Euclidean metric $dx^2+dy^2$ by a lattice
$e^{-\alpha}\Gamma_{\nu}$. Thus replacing $\Gamma_{\nu}$ (resp. $T_{\nu}$)
by $e^{-\alpha}\Gamma_{\nu}$ (resp. $T_{\nu}^{\prime}$) for a
sufficiently large $\alpha$, we may
initially assume that $||\Gamma_{\nu}||<1$, or
equivalently $||\Gamma_{\nu}^{*}|| >1$. Taking $a$ sufficiently large we
also assume that $e^{2a}>32$. Let $\omega=fdx+gdy+hdu$ be an element of
$C^{\infty}_c(Y_{a,\nu},\Omega^1(\chi_{\nu,i}))$. Then a computation in
{\bf Proposition 3.1} implies
\begin{equation}
\int_{Y_{a,\nu}}\Delta_Tf\cdot\bar{f}dxdydu \geq
 ||\Gamma_{\nu}^{*}||^2\int_{Y_{a,\nu}}|f|^2dxdydu \geq \int_{Y_{a,\nu}}|f|^2dxdydu,
\end{equation}
\begin{equation}
\int_{Y_{a,\nu}}\Delta_Tg\cdot\bar{g}dxdydu \geq
 ||\Gamma_{\nu}^{*}||^2\int_{Y_{a,\nu}}|g|^2dxdydu \geq \int_{Y_{a,\nu}}|g|^2dxdydu,
\end{equation}
and
\begin{equation}
\int_{Y_{a,\nu}}\Delta_Th\cdot\bar{h}e^{-2u}dxdydu \geq
 ||\Gamma_{\nu}^{*}||^2\int_{Y_{a,\nu}}|h|^2e^{-2u}dxdydu.
\end{equation}
Using
\[
 ||dx||=||dy||=e^u, \quad ||du||=1
\]
and {\bf Lemma 3.2}, an integration by parts shows
\begin{eqnarray*}
\int_{Y_{a,\nu}}(\Delta^1\omega,\omega)dv_g &=&
 \int_{Y_{a,\nu}}e^{2u}(\Delta_Tf\cdot\bar{f}+\Delta_Tg\cdot\bar{g})dxdydu\\
&+ &\int_{Y_{a,\nu}}|\nabla_Th|^2dxdydu\\
&+&
 \int_{Y_{a,\nu}}(|\partial_uf|^2+|\partial_ug|^2+|\partial_uh|^2e^{-2u})dxdydu\\
&+& 2\int_{Y_{a,\nu}}\{(\partial_xh\cdot\bar{f}+\partial_x\bar{h}\cdot f)+(\partial_yh\cdot\bar{g}+\partial_y\bar{h}\cdot g)\}dxdydu
\end{eqnarray*} 
\begin{eqnarray}
&=&
 \int_{Y_{a,\nu}}(e^{2u}-16)(\Delta_Tf\cdot\bar{f}+\Delta_Tg\cdot\bar{g})dxdydu\\
&+&
 16\int_{Y_{a,\nu}}\{(\Delta_Tf\cdot\bar{f}-|f|^2)+(\Delta_Tg\cdot\bar{g}-|g|^2)\}dxdydu\\
&+&
 \frac{1}{4}\int_{Y_{a,\nu}}\{64|f|^2+8(\partial_xh\cdot\bar{f}+\partial_x\bar{h}\cdot
 f)+|\nabla_Th|^2\}dxdydu\\
&+&
 \frac{1}{4}\int_{Y_{a,\nu}}\{64|g|^2+8(\partial_yh\cdot\bar{g}+\partial_y\bar{h}\cdot
 g)+|\nabla_Th|^2\}dxdydu\\
&+& \frac{1}{2}\int_{Y_{a,\nu}}|\nabla_Th|^2dxdydu\\
&+& \int_{Y_{a,\nu}}(|\partial_uf|^2+|\partial_ug|^2+|\partial_uh|^2e^{-2u})dxdydu.
\end{eqnarray}
By (15) and (16), (19) is nonnegative. Moreover  
\[
 (8|f|-|\nabla_Th|)^2 \leq 64|f|^2+8(\partial_xh\cdot\bar{f}+\partial_x\bar{h}\cdot
 f)+|\nabla_Th|^2
\] 
and
\[
 (8|g|-|\nabla_Th|)^2 \leq 64|g|^2+8(\partial_yh\cdot\bar{g}+\partial_y\bar{h}\cdot
 g)+|\nabla_Th|^2
\]
imply that both (20) and (21) are nonnegative. Since
\begin{eqnarray*}
 \int_{Y_{a,\nu}}(e^{2u}-16)(\Delta_Tf\cdot\bar{f}+\Delta_Tg\cdot\bar{g})dxdydu
  &=& \int_{Y_{a,\nu}}(e^{2u}-16)(|\nabla_Tf|^2+|\nabla_Tg|^2)dxdydu\\
& \geq& (e^{2a}-16)\int_{Y_{a,\nu}}(|\nabla_Tf|^2+|\nabla_Tg|^2)dxdydu\\
&=& (e^{2a}-16)\int_{Y_{a,\nu}}(\Delta_Tf\cdot\bar{f}+\Delta_Tg\cdot\bar{g})dxdydu
\end{eqnarray*}
we obtain
\begin{eqnarray*}
\int_{Y_{a,\nu}}(\Delta^1\omega,\omega)dv_g &\geq &
 (e^{2a}-16)\int_{Y_{a,\nu}}(\Delta_Tf\cdot\bar{f}+\Delta_Tg\cdot\bar{g})dxdydu\\
&+&  \frac{1}{2}e^{2a}\int_{Y_{a,\nu}}|\nabla_Th|^2e^{-2u}dxdydu\\
&\geq & \frac{1}{2}e^{2a}\int_{Y_{a,\nu}}(\Delta_Tf\cdot\bar{f}+\Delta_Tg\cdot\bar{g}+\Delta_Th\cdot\bar{h}e^{-2u})dxdydu.
\end{eqnarray*}
Here we have used the fact $e^{2a}$ is greater than $32$. Thus by (15), (16)
and (17) we see
\[
 \int_{Y_{a,\nu}}(\Delta^1\omega,\omega)dv_g \geq \frac{1}{2}e^{2a}||\Gamma_{\nu}^{*}||^2\int_{Y_{a,\nu}}||\omega||^2dv_g.
\]
Now {\bf Fact 3.1} yields
\begin{prop} For a sufficiently large $a$,
\[
 \mu_{1}(\Delta_{dir}^1|_{Y_{a,\nu}}) \geq \frac{1}{2}e^{2a}||\Gamma_{\nu}^{*}||^2.
\]
\end{prop}

By the same argument we will find
\begin{prop}Suppose $\alpha >1$. Then for a sufficiently large $a$ and
 $*=abs$ or $rel$, 
\[
 \mu_{1}(\Delta_{dir/*}^1|_{X_{\alpha a}\cap Y_{a, \nu}}) \geq \frac{1}{2}e^{2a}||\Gamma_{\nu}^{*}||^2.
\]
\end{prop}

\subsection{A convergence of spectrum}

In {\bf Lemma 3.5} we have shown that $\mu_n(\Delta^p_{dir}|_{X_a})$ is a
monotone decreasing function of $a$ which is bounded below by
$\mu_n(\Delta^p)$. In this section we will show the following theorems.
\begin{thm}
\[
 \lim_{a\to\infty}\mu_n(\Delta^p_{dir}|_{X_a})=\mu_n(\Delta^p).
\]
\end{thm}
\begin{thm}
\[
 \lim_{a\to\infty}\mu_n(\Delta^p_{*}|_{X_a})=\mu_n(\Delta^p)\quad
 *=abs, rel.
\]
\end{thm}

\begin{cor} Let $t$ be a positive number. Then 
\[
 {\rm Tr}[e^{-t\Delta^p}]=\lim_{a\to\infty}{\rm
 Tr}[e^{-t\Delta^p_{dir}|_{X_a}}]=\lim_{a\to\infty}{\rm
 Tr}[e^{-t\Delta^p_{*}|_{X_a}}],\quad *=abs, rel.
\]
\end{cor}
Let us fix a positive $a_0$ so that $Y_{a_0}$ is a disjoint union:
\[
 Y_{a_0}=\amalg_{\nu=1}^hT_{\nu}\times[a_0, \infty).
\]
We may assume that $e^{2a_0}>32$ and $||\Gamma_{\nu}^*||>1$ for
every $\nu$. Thus {\bf Proposition 3.1} and {\bf Proposition 3.3} are
available for $a>a_0$. Let us fix such an $a$ and let $\chi$ be a smooth function on $X$ satisfying
\begin{enumerate}
\item $0 \leq \chi \leq 1.$
\item $\chi|_{X_a}=1$ and $\chi|_{Y_{2a}}=0$.
\item $|\nabla \chi| \leq a^{-1}.$
\end{enumerate}
Let $\varphi_i$ be its eigenform of $\Delta^p$ whose eigenvalue is
$\mu_i(\Delta^p)$ and $\frak{M}_n$ an element of ${\rm
Gr}_nD(\Delta^p)$ spanned by $\{\varphi_1,\cdots,\varphi_n\}$. Then
for an arbitrary $\varphi\in \frak{M}_n$ we have
\begin{equation}
\int_X ||\nabla\varphi||^2dv_g = \int_X (\Delta^{p}\varphi, \varphi)dv_g
 \leq \mu_n(\Delta^p)\int_X ||\varphi||^2dv_g.
\end{equation}
The LHS is
\begin{eqnarray*}
\int_X ||\nabla\varphi||^2dv_g &=& \int_X
 ||\nabla(\chi\varphi)+\nabla((1-\chi)\varphi)||^2dv_g \\
&=& \int_X||\nabla(\chi\varphi)||^2dv_g +
 \int_X||\nabla((1-\chi)\varphi)||^2dv_g\\
&+& 2{\rm Re}\int_X(\nabla(\chi\varphi), \nabla((1-\chi)\varphi))dv_g.
\end{eqnarray*}
Since
\[
 \chi(1-\chi)\leq \frac{1}{4}\quad\mbox{and}\quad |\nabla\chi|\leq \frac{1}{a}
\]
and by Schwartz inequality,
\[
 |(\nabla(\chi\varphi), \nabla((1-\chi)\varphi))| \leq (\frac{1}{a}+\frac{1}{a^2})||\varphi||^2+(\frac{1}{a}+\frac{1}{4})||\nabla\varphi||^2.
\]
Therefore (24) implies
\begin{eqnarray*}
\mu_n(\Delta^p)\int_X ||\varphi||^2dv_g &\geq&
 \int_X||\nabla((1-\chi)\varphi)||^2dv_g+\int_X||\nabla(\chi\varphi)||^2dv_g\\
&-&2(\frac{1}{a}+\frac{1}{a^2})\int_X||\varphi||^2dv_g
 -2(\frac{1}{a}+\frac{1}{4})\int_X||\nabla\varphi||^2dv_g\\
& \geq & \int_X||\nabla((1-\chi)\varphi)||^2dv_g\\
&-&2\{\frac{1}{a}+\frac{1}{a^2}+\mu_n(\Delta^p)(\frac{1}{a}+\frac{1}{4})\}\int_X||\varphi||^2dv_g.
\end{eqnarray*}
Notice that $(1-\chi)\varphi$ is contained in the domain of
$\Delta^p_{dir}|_{Y_a}$. By {\bf Fact 3.2} and {\bf Proposition
3.1}(if $p=0$), or {\bf Propostion 3.3}(if $p=1$),
\begin{eqnarray*}
\int_X||\nabla((1-\chi)\varphi)||^2dv_g &\geq&
 \mu_1(\Delta^p_{dir}|_{Y_{a}})\int_X||(1-\chi)\varphi||^2dv_g\\
&\geq& \mu_1(\Delta^p_{dir}|_{Y_{a}})\int_{Y_{2a}}||\varphi||^2dv_g\\
&\geq&Ce^{2a}\int_{Y_{2a}}||\varphi||^2dv_g,
\end{eqnarray*}
where $C$ is a positive constant independent of $a$. So we have obtained
\[
 \{\mu_n(\Delta^p)+2(\frac{1}{a}+\frac{1}{a^2}+\mu_n(\Delta^p)(\frac{1}{a}+\frac{1}{4}))\}\int_X||\varphi||^2dv_g\\
 \geq Ce^{2a}\int_{Y_{2a}}||\varphi||^2dv_g.
\]
Now put
\[
 \rho_n(a)=2C^{-1}e^{-2a}\{(\frac{1}{a}+\frac{1}{a^2})+\mu_n(\Delta^p)(\frac{1}{a}+\frac{3}{4})\},
\]
then we have proved the following proposition.
\begin{prop}
For $\varphi\in\frak{M}_n$
\[
 \int_{Y_{2a}}||\varphi||^2dv_g \leq \rho_n(a)\int_X||\varphi||^2dv_g.
\]
\end{prop}
Let $\alpha$ be greater
than two and $\phi_i$ an eigenform of $\Delta^p_{*}|_{X_{\alpha
a}}$($*=rel, \mbox{or} \,abs$) whose eigenvalue is $\mu_i(\Delta^p_{*}|_{X_{\alpha
a}})$ and $\frak{M}_n(\alpha a)$ an element of ${\rm
Gr}(\Delta^p_{*}|_{X_{\alpha a}})$ spanned by
$\{\phi_1,\cdots,\phi_n\}$. 

\begin{prop} For $\phi\in \frak{M}_n(\alpha a)$,
\[
 \int_{Y_{2a\cap X_{\alpha a}}}||\phi||^2dv_g \leq
 \rho^{0}_n(a)\int_{X_{\alpha a}}||\phi||^2dv_g,
\]
where
\[
 \rho^{0}_n(a)=2C^{-1}e^{-2a}\{(\frac{1}{a}+\frac{1}{a^2})+\mu_n(\Delta^p_{dir}|_{X_{a_0}})(\frac{1}{a}+\frac{3}{4})\}.
\]
\end{prop}
{\bf Proof.} The argument is almost same as one of {\bf Proposition 3.5}. For $\phi\in \frak{M}_n(\alpha a)$,
\begin{equation}
\int_{X_{\alpha a}}||\nabla\phi||^2dv_g=\int_{X_{\alpha
 a}}(\Delta^p\phi,\phi)dv_g \leq \mu_n(\Delta^p_{*}|_{X_{\alpha a}})\int_{X_{\alpha a}}||\phi||^2dv_g.
\end{equation}
Using {\bf Proposition 3.2} and {\bf Proposition 3.4} instead {\bf
Proposition 3.1} and {\bf Proposition 3.3}, respectively the previous
computation will show
\begin{eqnarray*}
\mu_n(\Delta^p_{*}|_{X_{\alpha a}})\int_{X_{\alpha
 a}}||\phi||^2dv_g &\geq&  \int_{X_{\alpha a}}||\nabla\phi||^2dv_g\\
 &\geq& Ce^{2a}\int_{Y_{2a}\cap
 X_{\alpha a}}||\phi||^2dv_g\\
&-& 2\{\frac{1}{a}+\frac{1}{a^2}+\mu_n(\Delta^p_{*}|_{X_{\alpha
 a}})(\frac{1}{a}+\frac{1}{4})\}\int_{X_{\alpha a}}||\phi||^2dv_g,
\end{eqnarray*}
which yields
\[
 \int_{Y_{2a}\cap X_{\alpha a}}||\phi||^2dv_g \leq 2C^{-1}e^{-2a}\{(\frac{1}{a}+\frac{1}{a^2})+\mu_n(\Delta^p_{*}|_{X_{\alpha
 a}})(\frac{1}{a}+\frac{3}{4})\}\int_{X_{\alpha a}}||\phi||^2dv_g.
\]
By {\bf Lemma 3.5},
\[
 \mu_n(\Delta^p_{*}|_{X_{\alpha a}}) \leq \mu_n(\Delta^p_{dir}|_{X_{a_0}})
\]
and we have proved the proposition.
\begin{flushright}
$\Box$
\end{flushright}

{\bf A proof of Theorem 3.1}\\

As before let $\varphi_i$ be an eigenvector of $\Delta^p$ whose
eigenvalue is $\mu_i(\Delta^p)$ and $\frak{M}_{n,\chi}$ an
$n$-dimensional subspace of $D(\Delta^p_{dir}|_{X_{2a}})$ spanned by
$\{\chi\varphi_1,\cdots,\chi\varphi_n\}$. Let us choose
$\varphi\in\frak{M}_n$ to be
\[
 \frac{\int_X||\nabla(\chi\varphi)||^2dv_g}{\int_X||\chi\varphi||^2dv_g}=\sup_{f\in\frak{M}_{n,\chi}}\frac{\int_X||\nabla
 f||^2dv_g}{\int_X||f||^2dv_g}.
\]
By {\bf Lemma3.3} the RHS is greater than or equal to
$\mu_n(\Delta^p_{dir}|_{X_{2a}})$ and therefore
\[
 \int_X||\nabla(\chi\varphi)||^2dv_g \geq \mu_n(\Delta^p_{dir}|_{X_{2a}})\int_X||\chi\varphi||^2dv_g.
\]
On the other hand by a choice of $\chi$,
\begin{eqnarray*}
||\nabla(\chi\varphi)||^2 &\leq& ||\nabla\chi\cdot\varphi||^2+2|{\rm
 Re}(\nabla\chi\cdot\varphi,\chi\nabla\varphi)|+||\chi\nabla\varphi||^2\\
&\leq& \frac{1}{a^2}||\varphi||^2+\frac{2}{a}|{\rm Re}(\varphi,\nabla\varphi)|+||\nabla\varphi||^2\\
&\leq& (\frac{1}{a^2}+\frac{1}{a})||\varphi||^2+(\frac{1}{a}+1)||\nabla\varphi||^2
\end{eqnarray*}
hence
\begin{eqnarray*}
& &(\frac{1}{a^2}+\frac{1}{a})\int_X||\varphi||^2dv_g
 + (\frac{1}{a}+1)\int_X||\nabla\varphi||^2dv_g\\
&\geq& \mu_n(\Delta^p_{dir}|_{X_{2a}}) \int_X||\chi\varphi||^2dv_g\\
&\geq&  \mu_n(\Delta^p_{dir}|_{X_{2a}}) \int_{X_{a}}||\varphi||^2dv_g\\
&=& \mu_n(\Delta^p_{dir}|_{X_{2a}})(\int_{X}||\varphi||^2dv_g-\int_{Y_{a}}||\varphi||^2dv_g).
\end{eqnarray*}
Here by (24) the most left hand side is less than or equal to 
\[
 \{(\frac{1}{a}+1)\mu_n(\Delta^p)+(\frac{1}{a^2}+\frac{1}{a})\}\int_X||\varphi||^2dv_g,
\]
and by {\bf Proposition 3.5} the most right hand side is greater than or equal to
\[
 \mu_n(\Delta^p_{dir}|_{X_{2a}})(1-\rho_n(\frac{a}{2}))\int_X||\varphi||^2dv_g.
\]
Thus we have obtained
\[
 (\frac{1}{a}+1)\mu_n(\Delta^p)\geq \mu_n(\Delta^p_{dir}|_{X_{2a}})(1-\rho_n(\frac{a}{2}))-(\frac{1}{a^2}+\frac{1}{a}).
\]
Now since
\[
 \lim_{a\to\infty}\rho_n(\frac{a}{2})=0,
\]
and {\bf Lemma 3.5} implies the desired result.
\begin{flushright}
$\Box$
\end{flushright}

{\bf A proof of Theorem 3.2.}\\

 Since a proof is almost same as one of {\bf
Theorem 3.1} we will only indicate where a modification is necessary. As
before let $\phi_i$ be an eigenform of $\Delta^p_{*}|_{X_{3a}}$
whose eigenvalue is $\mu_i(\Delta^p_{*}|_{X_{3a}})$ and
$\frak{M}_n(3a)_{\chi}$ a $n$-dimensional subspace of
$D(\Delta^p_{dir}|_{X_{2a}})$ spanned by
$\{\chi\phi_1,\cdots,\chi\phi_n\}$. We choose $\phi\in\frak{M}_n(3a)$ so
that
\[
 \frac{\int_{X_{3a}}||\nabla(\chi\phi)||^2dv_g}{\int_{X_{3a}}||\chi\phi||^2dv_g}=\sup_{f\in\frak{M}_{n}(3a)_{\chi}}\frac{\int_X||\nabla
 f||^2dv_g}{\int_X||f||^2dv_g}.
\]
Then {\bf Lemma 3.3} implies
\[
  \int_{X_{3a}}||\nabla(\chi\phi)||^2dv_g \geq
  \mu_n(\Delta^p_{dir}|_{X_{2a}})\int_{X_{3a}}||\chi\phi||^2dv_g.
\]
Using (25) and {\bf Proposition 3.6} instead (24) and {\bf Proposition
3.5}, respectively the same computation as in {\bf Theorem 3.1} will yield
\[
  (\frac{1}{a}+1)\mu_n(\Delta^p_{*}|_{X_{3a}})\geq \mu_n(\Delta^p_{dir}|_{X_{2a}})(1-\rho_n^0(\frac{a}{2}))-(\frac{1}{a^2}+\frac{1}{a}).
\]
By {\bf Lemma 3.5} $\mu_n(\Delta^p_{*}|_{X_{3a}})$ is bounded by
$\mu_n(\Delta^p_{dir}|_{X_{3a}})$ from above.  Therefore
\[
 \lim_{a\to\infty}\mu_n(\Delta^p_{*}|_{X_{3a}})=\lim_{a\to\infty}\mu_n(\Delta^p_{dir}|_{X_{3a}}),
\]
and this is equal to $\mu_n(\Delta^p)$ by {\bf Theorem 3.1}.
\begin{flushright}
$\Box$
\end{flushright}

\subsection{A theorem of Cheeger-M\"{u}ller type}

Since both $\Delta^p$ and $\Delta^p_{abs}|_{X_{a}}$ have only discrete
spectrum the $p$-th cohomology groups $H^{p}(X_a,\rho)$ and
$H^{p}(X,\rho)$ is isomorphic to ${\rm Ker}\,\Delta^{p}_{abs}|_{X_a}$
and ${\rm Ker}\,\Delta^{p}_{X}$ by Hodge theory, respectively. If $a$ is sufficiently
large, $H^{p}(X_a,\rho)$ is isomorphic to $H^{p}(X,\rho)$ by the restriction
and therefore ${\rm Ker}\,\Delta^{p}_{X}$ is also isomorphic to
${\rm Ker}\,\Delta^{p}_{abs}|_{X_a}$. 
Since $\rho$ is cuspidal $H^0(T_{\nu},\rho)=0$ for every $\nu$ and 
\[
 H^2(T_{\nu},\rho)=0,
\]
by Poincar\'{e} duality. Because $\rho$ is a local system on a flat
torus $T_{\nu}$ we see 
\[
 H^1(T_{\nu},\rho)=0,
\]
by the index theorem. Therefore the exact sequence
\begin{eqnarray*}
& & H^{p-1}(\partial X_a,\rho)=\oplus_{\nu=1}^h H^{p-1}(T_{\nu},\rho) \to
 H^p(X_a,\partial X_a, \rho)\to H^p(X_a, \rho)\\
 &\to& H^{p}(\partial X_a,\rho)=\oplus_{\nu=1}^h H^{p}(T_{\nu},\rho)
\end{eqnarray*} 
shows $H^p(X_a,\partial X_a, \rho)$ is isomorphic to $H^p(X_a,
\rho)$. Thus we find
\[
 h^{p}(X,\rho)=h^{3-p}(X,\rho),
\]
by Poincar\'{e} duality. (Recall $h^{p}(X,\rho)$ is the dimension of $H^p(X_a,
\rho)$.) In particular since $h^{0}(X,\rho)$ vanishes so does $h^{3}(X,\rho)$. 
Moreover Hodge $*$ operator yields an
isomorphism
\[
 {\rm Ker}\,\Delta^{p}_{abs}|_{X_a}\stackrel{*}\simeq {\rm Ker}\,\Delta^{3-p}_{rel}|_{X_a},
\]
and therefore
\[
 h^{p}(X,\rho)={\rm dim}{\rm Ker}\,\Delta^{p}_{abs}|_{X_a}={\rm Ker}\,\Delta^{p}_{rel}|_{X_a}.
\]
{\it A partial spectral zeta function} of $\Delta_{*}^{p}|_{X_a}$
($*=abs$ or $rel$) is defined to be
\[
 \zeta_{X_a,*}^{(p)}(z,\rho)=\frac{1}{\Gamma(z)}\int^{\infty}_{0}\theta_p(t,a)t^{z-1}dt,\quad \theta_p(t,a)={\rm
 Tr}[e^{-t\Delta_{*}^{p}|_{X_a}}]-h^{p}(X,\rho).
\]
(Here $a$ is assumed to be sufficiently large.) If ${\rm Re}\,z$ is sufficiently large it absolutely
converges. Since $\theta_p(t,a)$
has an asymptotic expansion
\begin{equation}
  \theta_p(t,a)\sim
  t^{-\frac{3}{2}}\sum_{l=0}^{2N}c_l(a)t^{l/2}+O(t^{N-\frac{3}{2}}), \quad
  \mbox{as}\quad t\to 0.
\end{equation}
the same argument as the beginning of this section will show that
$\zeta_{X_a,*}^{(p)}(z,\rho)$ is meromorphically continued to the whole
plane and that it is regular at the
origin. We define {\it the spectral zeta function} of $X_a$ as
\[
  \zeta_{X_a}(z,\rho)=\sum_{p=0}^{3}(-1)^p p\cdot \zeta_{X_a,abs}^{(p)}(z,\rho).
\]

Since Hodge $*$ operator commutes with Hodge Laplacian and since it interchanges two boundary conditions,
\[
 \zeta_{X}^{(p)}(z,\rho)=\zeta_{X}^{(3-p)}(z,\rho),\quad \zeta_{X_a,abs}^{(p)}(z,\rho)=\zeta_{X_a,rel}^{(3-p)}(z,\rho).
\]

Therefore
\[
 \zeta_{X_a}(z,\rho)=2\zeta_{X_a,rel}^{(1)}(z,\rho)-\zeta_{X_a,abs}^{(1)}(z,\rho)-3\zeta_{X_a,rel}^{(0)}(z,\rho)
\]
and
\[
 \zeta_{X}(z,\rho)=\zeta_{X}^{(1)}(z,\rho)-3\zeta_{X}^{(0)}(z,\rho).
\]
\begin{thm}
\[
 \lim_{a\to \infty} \zeta_{X_a}(z,\rho)=\zeta_{X}(z,\rho).
\]
\end{thm}

\begin{cor}
\[
 \lim_{a\to \infty}\zeta_{X_a}^{\prime}(0,\rho)=\zeta_{X}^{\prime}(0,\rho).
\]
\end{cor}

{\bf Proof of Theorem 3.3.}
Let us write 
\[
 \int^{\infty}_0 \theta_p(t)t^{z-1}dt=\int^{1}_0 \theta_p(t)t^{z-1}dt+\int^{\infty}_1 \theta_p(t)t^{z-1}dt,
\]
\[
 \int^{\infty}_0 \theta_p(t,a)t^{z-1}dt=\int^{1}_0 \theta_p(t,a)t^{z-1}dt+\int^{\infty}_1 \theta_p(t,a)t^{z-1}dt
\]
We will investigate convergence of corresponding terms in RHS. For the second
term {\bf Theorem 3.2} implies
\[
 \lim_{a\to \infty}\int^{\infty}_1 \theta_p(t,a)t^{z-1}dt=\int^{\infty}_1 \theta_p(t)t^{z-1}dt.
\]
As we have seen at the begining of this section $\theta_p(t)$ has an
asymptotic expansion (see (14)),
\[
 \theta_p(t)\sim
 t^{-\frac{3}{2}}\sum_{l=0}^{2N}c_lt^{l/2}+O(t^{N-\frac{3}{2}}),\quad c_{2l+1}~=0,
\] 
around $t=0$. We  put
\[
 \theta_p(t,a)^{(N)}=\theta_p(t,a)-t^{-\frac{3}{2}}\sum_{l=0}^{2N}c_l(a)t^{l/2}, \quad \theta_p(t)^{(N)}=\theta_p(t)-t^{-\frac{3}{2}}\sum_{l=0}^{2N}c_lt^{l/2}.
\]
Then $t^{\frac{3}{2}-N}\theta_p(t,a)^{(N)}$ and
$t^{\frac{3}{2}-N}\theta_p(t)^{(N)}$ are bounded on $(0,1]$. 
By {\bf Theorem 3.2} we see 
\[
 \lim_{a\to\infty}c_l(a)=c_l,\quad 0\leq l\leq 2N,
\]
which implies in turn
\[
 \lim_{a\to\infty}\theta_p(t,a)^{(N)}=\theta_p(t)^{(N)}\quad
 \mbox{on}\quad (0,1]
\]
Therefore if ${\rm Re}z >\frac{3}{2}-N$,
\[
 \lim_{a\to \infty}\int^{1}_0\theta_p(t,a)^{(N)}t^{z-1}dt=\int^{1}_0\theta_p(t)^{(N)}t^{z-1}dt.
\]
and 
\[
 \int^{1}_0 \theta_p(t,a)t^{z-1}dt=\sum_{l=0}^{2N}\frac{c_l(a)}{z+\frac{l-3}{2}}+\int^{1}_0\theta_p(t,a)^{(N)}t^{z-1}dt,
\]
converges to 
\[
 \int^{1}_0 \theta_p(t)t^{z-1}dt=\sum_{l=0}^{2N}\frac{c_l}{z+\frac{l-3}{2}}+\int^{1}_0\theta_p(t)^{(N)}t^{z-1}dt,
\]
as $a\to\infty$.

\begin{flushright}
$\Box$
\end{flushright}

For a finite dimensional vector space $V$ we set
\[
 {\rm det}V=\wedge^{{\rm dim} V}V.
\]
{\it The determinant} of a bounded complex of finite dimensional
vector spaces $(C^{\cdot},\,\partial)$ is defined to be
\[
 {\rm det}(C^{\cdot},\,\partial)=\otimes_{i} ({\rm det}C^{i})^{(-1)^{i}}.
\] 
Here for a one dimensional complex vector space $L$, $L^{-1}$ is its
dual. Due to Knudsen and Mumford, there is a canonical isomorphism
\begin{equation}
 {\rm det}(C^{\cdot},\,\partial)\simeq \otimes_{i}{\rm det}\,H^{i}(C^{\cdot},\,\partial)^{(-1)^i}.
\end{equation}
Let $\Sigma=\{\Sigma_p\}_{p}$ be a triangulation of $X_a$ where
$\Sigma_p$ is the set of $p$-simplices and ${\bf e}=\{{\bf
e}_1,\cdots,{\bf e}_r\}$ a unitary basis of $\rho$. We define a Hermitian
inner product on the group of $p$-cochains:
\[
 C^{p}(\Sigma,\rho)=C^{p}(\Sigma)\otimes \rho,
\]
so that $\{[\sigma]^{*}\otimes{\bf e}_i\}$ form its unitary basis, where
$[\sigma]^{*}$ is the dual vector of $[\sigma]$. Now the Knudsen and
Mumford isomorphism induces a
metric $||\cdot||_{FR,a}$ on ${\rm
det}H^{\cdot}(X_a,\rho)=\otimes_{i}{\rm
det}\,H^{i}(C^{\cdot}(\Sigma,\rho))^{(-1)^i}$, which is referred as {\it
Franz-Reidemeister metric}. Via the isomorphism
\[
 H^{\cdot}(X_a,\rho)\simeq H^{\cdot}(X,\rho),
\]
it yields a metric $||\cdot||_{FR,a}$ on ${\rm det}H^{\cdot}(X,\rho)$. Notice that they are
independent of $a$ as far as it is sufficiently large since we can
use the same triangulations to define them. Therefore the limit
\[
 ||\cdot||_{FR}=\lim_{a\to \infty}||\cdot||_{FR,a}
\]
is well-defined. For a later purpose we will describe it in terms of a
combinatric zeta function.\\

 A triangulation $\Sigma$ of $X_a$ induces a simplicial decomposition $\tilde{\Sigma}$ of
the universal covering $\tilde{X_a}$. Let $\{\sigma^{(p)}_1,\cdots,\sigma^{(p)}_{\gamma_{p}}\}$ the
set of $p$-simplices of $\tilde{\Sigma}$ which are a lift of $\Sigma_p$. Then $C_{p}(\tilde{\Sigma})$ is a free ${\mathbb
C}[\Gamma]$-module genereted by these elements. A twisted chain complex
is defined to be
\[
 C_{\cdot}(\Sigma,\rho)=C_{\cdot}(\tilde{\Sigma})\otimes_{{\mathbb C}[\Gamma]}\rho,
\]
which is a bounded complex of finite dimensional vector spaces. We will
introduce a Hermitian inner product so that $\{\sigma^{(p)}_i\otimes{\bf
e}_j\}$ is a unitary basis. Here is an explict description of the
boundary map: Let $\sum_{k}(-1)^k\gamma_k[\sigma^{(p-1)}_{k}]\, (\gamma_k\in \Gamma)$ be the boundary of
$[\sigma^{(p)}_i]\in C^{p}(\tilde{\Sigma})$. Then 
\[
 \partial([\sigma^{(p)}_i]\otimes{\bf e}_j)=\sum_{k}(-1)^k
 [\sigma^{(p-1)}_{k}]\otimes\rho(\gamma_k){\bf e}_j.
\]
Let $(C^{\cdot}(\Sigma,\rho),\,\delta)$ be the dual complex. By the
inner product we may identify $C^{\cdot}(\Sigma,\rho)$ with
$C_{\cdot}(\Sigma,\rho)$ and in particular the dual vector of
$[\sigma^{(p)}_i]\otimes{\bf e}_j$ will be identified with itself. Thus
$(C^{\cdot}(\Sigma,\rho),\,\delta)$ is a complex such that each
$C^{p}(\Sigma,\rho)$ is nothing but
$C_{p}(\Sigma,\rho)$ as a vector space but the differential $\delta$ is the
Hermitian dual of $\partial.$ Let us define {\it a (positive) combinatric
Laplacian} $\Delta^p_{comb}$ on $C^p(\Sigma, \rho)(=C_p(\Sigma, \rho))$ to be
\[
 \Delta^p_{comb}=\partial \delta+\delta \partial.
\]
Then 
\[
 H^{p}(X_a,\rho)(=H_{p}(X_a,\rho))={\rm Ker}[\Delta_{comb}^p],
\]
has the inner product $(\cdot,\cdot)_{l^2,X_a}$ induced
by the metric on $C_{p}(\Sigma,\rho)$ (Here we have identified $H^{p}(X_a,\rho)$ with ${\rm
Ker}[\Delta_{comb}^p]$ which is a subspace of $C_{p}(\Sigma,\rho)$). It induces a norm
$|\cdot|_{l^2,X_a}$ on the determinant $\otimes_{p}{\rm
det}H^{p}(X_a,\rho)^{(-1)^p}$. Let us define {\it the
combinatric zeta function} to be
\[
  \zeta_{comb}(s,X_a)=\sum_{p}(-1)^p p\cdot \zeta^{(p)}_{comb}(s,X_a),
\]
where
\[
 \zeta^{(p)}_{comb}(s,X_a)=\sum_{\lambda}\lambda^{-s}.
\]					       
Here $\lambda$ runs through positive eigenvalues of $\Delta^p_{comb}$ on
$C^{p}(\Sigma,\rho)$.  {\it The modified Franz-Reidemeister torsion}
$\tau^{*}(X_a,\rho)$ is defined as 
\[
 \tau^{*}(X_a,\rho)=\exp (-\frac{1}{2}\zeta^{\prime}_{comb}(0,\rho)).
\]
If $H^{1}(X,\rho)$ vanishes so does every
$H^{p}(X,\rho)$ and $\tau^{*}(X_a,\rho)$ is the usual
Franz-Reidemeister torsion $\tau(X_a,\rho)$(\cite{Ray-Singer}). It is
known that $||\cdot||_{FR,a}$ is
equal to (\cite{Bismut-Zhang}\cite{Ray-Singer})
\[
 |\cdot|_{l^2,X_a}\cdot\tau^{*}(X_a,\rho).
\] 
By construction since both
$|\cdot|_{l^2,X_a}$ and $\tau^{*}(X_a,\rho)$ depend only on a
triangulation $\Sigma$, they are independent of sufficiently large $a$ as before. Thus the limit
\[
 |\cdot|_{l^2,X}=\lim_{a\to \infty}|\cdot|_{l^2,X_a}, \quad
  \tau^{*}(X,\rho)=\lim_{a\to \infty}\tau^{*}(X_a,\rho),
\] 
is well-defined and we set
\[
 ||\cdot||_{FR}=|\cdot|_{l^2,X}\cdot\tau^{*}(X,\rho).
\]

On the other hand since $H^{p}(X_a,\rho)$ is isomorphic to 
\[
 {\rm Ker} \Delta_{abs}^{p}|_{X_a}\subset L^{2}(X_a,\Omega^{p}(\rho))
\]
the inner product on $L^{2}(X_a,\Omega^{p}(\rho))$ induces a metric on
$H^{p}(X_a,\rho)$. Thus by the isomorphism $H^{p}(X,\rho)\simeq
H^{p}(X_a,\rho)$ we have a norm $|\cdot|_{L^2,X_a}$ on ${\rm det}H^{\cdot}(X,\rho)$.
{\it The Ray-Singer metric} $||\cdot||_{RS,a}$ is defined to be
\[
 ||\cdot||_{RS,a}=|\cdot|_{L^2,X_a}\cdot {\rm exp}(-\frac{1}{2}\zeta^{\prime}_{X_a}(0,\,\rho)).
\]
Similary using the canonical isomorphism
\[
 H^{p}(X,\rho)\simeq {\rm Ker} \Delta^{p}\subset L^{2}(X,\Omega^{p}(\rho)),
\]
{\it Ray-Singer metric} $||\cdot||_{RS}$ on ${\rm det}H^{\cdot}(X,\rho)$
is defined as
\[
 ||\cdot||_{RS}=|\cdot|_{L^2,X}\cdot {\rm exp}(-\frac{1}{2}\zeta^{\prime}_{X}(0,\,\rho)).
\]
Then we will show
\begin{thm}
$||\cdot||_{FR}$ and $||\cdot||_{RS}$ coincide. In particular
\[
 {\rm exp}(-\zeta^{\prime}_{X}(0,\rho))=\left(\frac{|\cdot|_{l^2,X}}{|\cdot|_{L^2,X}}\right)^2\tau^{*}(X,\rho)^2.
\]
\end{thm}

\begin{prop}
\[
 \lim_{a\to\infty}|\cdot|_{L^2,X_a}=|\cdot|_{L^2,X}.
\]
\end{prop}

In fact for a sufficiently large $a$ let
$\{\xi_{a,i}\}_{i}$ be an orthonormal base of ${\rm
Ker}\,\Delta_{abs}^{p}|_{X_a}$ and we define a map
\[
 {\rm Ker}\Delta^{p} \stackrel{P_a}\to {\rm Ker}\Delta_{abs}^{p}|_{X_a}
\]
as
\[
 P_a\psi=\sum_{i}\int_{X_a}(\psi,\xi_{a,i})dv_g\cdot \xi_{a,i}.
\]
Then we claim the following.
\begin{lm}
\[
 \lim_{a\to\infty}\int_{X_a}||\psi-P_a\psi||^2dv_g=0.
\]
\end{lm}
The following corollary will imply {\bf Proposition 3.7}.
\begin{cor} For $\psi\in {\rm Ker}\Delta^{p}$,
\[
 \lim_{a\to \infty}\int_{X_a}||P_a\psi||^2dv_g = \int_X||\psi||^2dv_g.
\]
\end{cor}
{\bf Proof.} Immediately from {\bf Proposition 3.5} and {\bf Lemma 3.6}.
\begin{flushright}
$\Box$
\end{flushright}

{\bf Proof of Lemma 3.6}. In the following arguments all $C$ are positive constants
independent of $a$. Let $\phi_{\lambda}$ be
an eigenform of $\Delta_{abs}^{p}|_{X_a}$ whose eigenvalue is
$\lambda$ and that  
\[
 \int_{X_a}||\phi_{\lambda}||^2dv_g=1.
\]
Let us expand $\psi$ as
\[
 \psi=\sum_{\lambda}\int_{X_a}(\psi,\phi_{\lambda})dv_g\cdot \phi_{\lambda}.
\]
Since 
\[
 \int_{X_a}||\psi-P_a\psi||^2dv_g = \sum_{\lambda > 0}|\int_{X_a}(\psi,\phi_{\lambda})dv_g|^2,
\]
it is sufficient to show that for $\phi=\phi_{\lambda}$ ($\lambda>0$), 
\[
 |\int_{X_a}(\psi,\phi)dv_g|\leq Ce^{-a}(\int_{X_a}||\psi||^2dv_g+C).
\]
Let us choose $\chi\in C^{\infty}_{c}(X_a)$ so that
\begin{enumerate}
\item $0\leq \chi \leq 1$.
\item $|\nabla \chi|$, $|\Delta\chi|$ are bounded by $1$.
\item $\chi|_{X_{a/2}}= 1$. 
\end{enumerate}
By Stokes theorem,
\begin{equation}
\int_{X_a}(\Delta^{p}(\chi\psi),\phi)dv_g =
 \int_{X_a}(\chi\psi,\Delta^{p}\phi)dv_g = \lambda \int_{X_a}\chi(\psi,\phi)dv_g
\end{equation}

Since $\Delta^{p}\psi=0$ and by the property 3 of $\chi$, LHS of (28) becomes
\begin{eqnarray*}
\int_{X_a}(\Delta^{p}(\chi\psi),\phi)dv_g &=&
 \int_{X_a}(\Delta\chi\cdot\psi,\phi)dv_g+2\int_{X_a}(\nabla\chi\cdot\nabla\psi,\phi)dv_g\\
&=& \int_{Y_{a/2}\cap X_a}(\Delta\chi\cdot\psi,\phi)dv_g+2\int_{Y_{a/2}\cap X_a}(\nabla\chi\cdot\nabla\psi,\phi)dv_g.
\end{eqnarray*}
Let us consider the first term. Using the property 2 of $\chi$ the
Schwartz inequality implies
\begin{eqnarray*}
 |\int_{Y_{a/2}\cap X_a}(\Delta\chi\cdot\psi,\phi)dv_g|&\leq&
  \frac{1}{2}(\int_{Y_{a/2}\cap X_a}||\psi||^2dv_g+\int_{Y_{a/2}\cap
  X_a}||\phi||^2dv_g)\\
&\leq& \frac{1}{2}(\int_{Y_{a/2}}||\psi||^2dv_g+\int_{Y_{a/2}\cap
  X_a}||\phi||^2dv_g).
\end{eqnarray*}
By {\bf Proposition 3.5},
\begin{eqnarray*}
\int_{Y_{a/2}}||\psi||^2dv_g &\leq& Ce^{-a/2} \int_{X}||\psi||^2dv_g\\
&=& Ce^{-a/2}(\int_{X_a}||\psi||^2dv_g+\int_{Y_a}||\psi||^2dv_g)\\
&\leq & Ce^{-a/2}\int_{X_a}||\psi||^2dv_g+Ce^{-a/2}\int_{Y_{a/2}}||\psi||^2dv_g,
\end{eqnarray*}
and therefore changing $C$ we obtain
\[
 \int_{Y_{a/2}}||\psi||^2dv_g \leq Ce^{-a/2}\int_{X_a}||\psi||^2dv_g.
\] 
Using {\bf Proposition 3.6} instead {\bf Proposition 3.5} the same
computation will show
\begin{equation}
 \int_{Y_{a/2}\cap X_a}||\phi||^2dv_g \leq Ce^{-a/2}\int_{X_a}||\phi||^2dv_g=Ce^{-a/2}
\end{equation}
and thus
\[
 |\int_{Y_{a/2}\cap X_a}(\Delta\chi\cdot\psi,\phi)dv_g| \leq Ce^{-a/2}(\int_{X_a}||\psi||^2dv_g+C).
\]

Next we will estimate the second term. Using the property 2 of $\chi$,
the Schwartz inequality yields,
\[
 |2\int_{X_a}(\nabla\chi\cdot\nabla\psi,\phi)dv_g|\leq
  \int_{Y_{a/2}}||\nabla\psi||^2dv_g+ \int_{Y_{a/2}\cap X_a}||\phi||^2dv_g.
\]
Since
\[
  \int_{Y_{a/2}}||\nabla\psi||^2dv_g \leq  \int_{X}||\nabla\psi||^2dv_g=\int_{X}(\psi,\Delta^{p}\psi)dv_g=0,
\]
and by (29), we obtain
\[
 |2\int_{X_a}(\nabla\chi\cdot\nabla\psi,\phi)dv_g|\leq Ce^{-a/2}.
\] 
Thus LHS of (28) is bounded by $Ce^{-a/2}(\int_{X_a}||\psi||^2dv_g+C)$.
Let us consider RHS of (28). The property 3 implies
\[
 \int_{X_a}\chi(\psi,\phi)dv_g=\int_{X_{a/2}}(\psi,\phi)dv_g+\int_{Y_{a/2}\cap X_a}\chi(\psi,\phi)dv_g.
\]
But by the previous arguments
\begin{eqnarray*}
|\int_{Y_{a/2}\cap X_a}\chi(\psi,\phi)dv_g| &\leq & \int_{Y_{a/2}\cap
 X_a}|(\psi,\phi)|dv_g\\
&\leq & \frac{1}{2}\{\int_{Y_{a/2}}||\psi||^2dv_g+\int_{Y_{a/2}\cap
 X_a}||\phi||^2dv_g\}\\
&\leq & Ce^{-a/2}(\int_{X_a}||\psi||^2dv_g+C),
\end{eqnarray*}
and we will
obtain
\[
 |\int_{X_{a/2}}(\psi,\phi)dv_g| \leq Ce^{-a/2}(\int_{X_{a}}||\psi||^2dv_g+C).
\]
Notice that
\begin{eqnarray*}
|\int_{X_a}(\psi,\phi)dv_g-\int_{X_{a/2}}(\psi,\phi)dv_g|&= &
 |\int_{X_a\cap Y_{a/2}}(\psi,\phi)dv_g|\\
 &\leq & \int_{X_a\cap Y_{a/2}}|(\psi,\phi)|dv_g\\
&\leq& Ce^{-a/2}(\int_{X_{a}}||\psi||^2dv_g+C),
\end{eqnarray*}
and the desired result has been obtained since
\begin{eqnarray*}
|\int_{X_a}(\psi,\phi)dv_g| &\leq &
 |\int_{X_a}(\psi,\phi)dv_g-\int_{X_{a/2}}(\psi,\phi)dv_g|+|\int_{X_{a/2}}(\psi,\phi)dv_g|\\
&\leq & Ce^{-a/2}(\int_{X_{a}}||\psi||^2dv_g+C).
\end{eqnarray*}

\begin{flushright}
$\Box$
\end{flushright}
Let us choose a sufficiently large number $a$ and small positive number $\delta$. Let $g_{0}$ be a Riemannian metric on $X$ such that
\[
 g_0(x)=
\left\{
\begin{array}{ccc}
g(x)&{\mbox if} & x\in X_{a-\delta}\\
du^2+e^{-2a}(dx^2+dy^2)&{\mbox if} & x\in Y_{a}
\end{array}
\right.
\]
We will consider a one parameter family of metrics:
\[
 g_q=(1-q)g_0+qg, \quad 0\leq q \leq 1.
\]
Let $\{{\bf e}^0,{\bf e}^1,{\bf e}^2\}$ be an orthonormal frame of
$\Omega^1|_{\partial X_a}$ with respect to $g(q)|_{\partial
X_a}=du^2+e^{-2a}(dx^2+dy^2)$ so that ${\bf e}^0=du$. Let $h(q)$ and $R(q)$
be the second fundamental of $\partial X_a$ with respect to $g(q)$ and the curvature tensor of
$g(q)$, respectively. Then we define elements

\[
 \hat{h}(q)=\sum_{1\leq a,b\leq 2}h(q)_{ab}{\bf
 e}^a\otimes{\bf e}^b
\]
and
\[
 \hat{R}_0(q)=\frac{1}{4}\sum_{j,k,l}R(q)_{0jkl}{\bf e}^j\otimes({\bf
 e}^k \wedge {\bf e}^l)
\]
of $\Omega^{\cdot}|_{\partial X_a}\otimes\Omega^{\cdot}|_{\partial
X_a}$. Using Berezin integral \cite{Bismut-Zhang}, $\int^{B}$, we put
\[
 \phi_a=\int^{1}_{0}dq\int^{B}\hat{h}(q)\hat{R}_0(q)\in
 \Omega^{\cdot}_{\partial X_a}.
\]
\begin{fact}(\cite{Dai-Fang})
\[
 \log
 \left(\frac{||\cdot||_{RS,a}}{||\cdot||_{FR,a}}\right)=\chi(\partial
 X_a,\rho)\log 2 +\gamma\cdot r \int_{\partial X_a}\phi_a,
\]
where $\gamma$ is an absolute constant.
\end{fact}
 Notice that the term
 $\tilde{e}(g_0,g_q)$ in the original formula vanishes because the
 dimension of $X$ is three. A straightforward computation will show
 that the norm of $\phi_a$ is bounded by a constant $C$ which is independent of $a$. Thus,
\[
 |\int_{\partial X_a}\phi_a| \leq C\cdot vol(\partial X_a)\leq
  C^{\prime}\cdot e^{-2a},
\]
where $C^{\prime}$ is also independent of $a$. 
Since $\rho$ is a unitary local system on $\partial X_a$ which is  a
disjoint union of flat tori, $\chi(\partial
X_a,\rho)$ vanishes by the index theorem. Therefore we have shown 

\begin{prop}
\[
 \lim_{a\to \infty}||\cdot||_{RS,a}=||\cdot||_{FR}.
\]
\end{prop}
{\bf Proof of Theorem 3.4.}
By {\bf Proposition 3.8}
\[
 ||\cdot||_{FR}=\lim_{a\to \infty}\{|\cdot|_{L^2,X_a}\cdot \exp(-\frac{1}{2}\zeta_{X_a}^{\prime}(0,\rho))\}.
\]
But by {\bf Propostion 3.7} and {\bf Corollary 3.2} this is equal to
$||\cdot||_{RS}$. 
\begin{flushright}
$\Box$
\end{flushright}

\subsection{A computation of the leading coefficient}
We will interprete the ratio $|\cdot|_{l^2,X}/|\cdot|_{L^2,X}$ as a
 period. Hereafter we will identify $H^p(X,\rho)$ and ${\rm
 Ker}\Delta_X^p$ by Hodge theory. Let $\phi^{(p)}=\{\phi^{(p)}_{1},\cdots,\phi^{(p)}_{h^p(X,\rho)}\}$ and
$\psi^{(p)}=\{\psi^{(p)}_{1},\cdots,\psi^{(p)}_{h^p(X,\rho)}\}$ be its
 unitary basis with respect to $(\,,\,)_{l^2,X}$ and $(\,,\,)_{L^2,X}$,
respectively. Then $\phi^{(p)}$ determines the dual basis
 $\phi_{(p)}=\{\phi_{(p),i}\}_{1\leq i \leq h^p(X,\rho)}$ of
 $H_p(X,\rho)$ and we write
\[
 \psi^{(p)}_{i}=\sum_{j=1}^{h^p(X,\rho)}\int_{\phi_{(p),j}}\psi^{(p)}_{i}\cdot \\\phi^{(p)}_{j}.
\]
{\it A period matrix of twisted $p$-forms} is defined to be
\[
 P(X,\rho)_p=(\int_{\phi_{(p),j}}\psi^{(p)}_{i})_{ij}.
\]
We call an alternating product $\prod_{p}|\det P(X,\rho)_p|^{(-1)^p}$ {\it a period of
$(X,\rho)$} and will denote it by ${\rm Per}(X,\rho)$. A simple
computation shows 
\[
 \psi^{(p)}_{1}\wedge \cdots \wedge \psi^{(p)}_{h^p(X,\rho)}=\det P(X,\rho)_p\cdot
 \phi^{(p)}_{1}\wedge \cdots \wedge \phi^{(p)}_{h^p(X,\rho)}.
\]
and by definition we have 
\[
 |\otimes_p(\psi^{(p)}_{1}\wedge \cdots \wedge \psi^{(p)}_{h^p(X,\rho)})^{(-1)^p}|_{L^2,X}=|\otimes_p(\phi^{(p)}_{1}\wedge \cdots \wedge \phi^{(p)}_{h^p(X,\rho)})^{(-1)^p}|_{l^2,X}=1
\]
Therefore  
\begin{eqnarray*}
\frac{|\otimes_p(\psi^{(p)}_{1}\wedge \cdots \wedge
 \psi^{(p)}_{h^p(X,\rho)})^{(-1)^p}|_{l^2,X}}{|\otimes_p(\psi^{(p)}_{1}\wedge
 \cdots \wedge \psi^{(p)}_{h^p(X,\rho)})^{(-1)^p}|_{L^2,X}} &=& |\otimes_p(\psi^{(p)}_{1}\wedge \cdots \wedge
 \psi^{(p)}_{h^p(X,\rho)})^{(-1)^p}|_{l^2,X}\\
&=& \prod_{p}|\det P(X,\rho)_p|^{(-1)^p}\\
&=& {\rm Per}(X,\rho). 
\end{eqnarray*}
Now using {\bf Fact 3.1}, {\bf Theorem 3.4} is reformulated as the following.

\begin{thm}
\[
 \lim_{z\to 0}z^{-2h^{1}(X,\rho)}R_{X}(z,\rho)=(\tau^{*}(X,\rho)\cdot {\rm Per}(X))^2.
\]
\end{thm}
\begin{cor} Suppose that $h^{1}(X,\rho)$ vanishes. Then
\[
 R_{X}(0,\rho)=\tau(X,\rho)^{2},
\]
where $\tau(X,\rho)$ is the usual Franz-Reidemeister torsion.
\end{cor}

\subsection{An example}

Let $K$ be a knot in $S^3$ whose complement $X_{K}$ admits a complete
hyperbolic structure of finite volume and $\rho$ a cuspidal unitary local system
of rank $r$ on $X_K$. Here the representaion of $\pi_1(X_K)$ associated
to $\rho$ is denoted by the same character. $X_K$ is obtained by
attaching 3-cells to a two dimensional CW-complex $L$ which is a
deformation retract of $X_K$. The argument of \cite{MilnorW}{\bf Lemma 7.2}
will show the following.
\begin{lm}
\[
 \tau(X_K,\rho)=\tau(L,\rho).
\]
\end{lm}
We will compute $\tau(L,\rho).$ Let
\[
 \pi_1(X_K)=<x_1,\cdots,x_n \,|\,r_1,\cdots, r_{n-1}>
\] 
be the Wirtinger presentation. Here $\{x_i\}_i$ (resp. $\{r_j\}_j$) is
generators (resp. relators). $H_1(X_K,{\mathbb Z})$ is an infinite cyclic
group and we fix a generator $t$. We choose $x_i$ so that Hurewicz map
\[
 \pi_1(X_K)\stackrel{\epsilon}\to H_1(X_K,{\mathbb Z})
\]
sends $x_i$ to $t$. Then a group ring ${\mathbb C}[H_1(X_K,{\mathbb Z})]$ is
isomorphic to Laurent polynomial ring $\Lambda={\mathbb C}[t,t^{-1}]$
and a ring homomorphism:
\[
 {\mathbb C}[\pi_1(X_K)]\to \Lambda.
\]
induced by Hurewicz map will be also denoted by $\epsilon$. 
Also the representation $\rho$ yields a homomorphism
\[
 {\mathbb C}[\pi_1(X_K)]\stackrel{\rho}\to M_{r}({\mathbb C})
\]
and let 
\[
 {\mathbb C}[\pi_1(X_K)]\stackrel{\epsilon\otimes\rho}\to M_{r}(\Lambda).
\]
be the tensor product of them. Composing this with the homomorphism
induced by the natural projection from
the free group $F_n$ of $n$-generators to $\pi_1(X_K)$ we obtain, 
\[
 {\mathbb C}[F_n]\stackrel{\Phi} \to M_{r}(\Lambda).
\]
The set of $0$-cells of $L$ consists of only one point $P_0$ and one of
1-cells is
\[
 \{x_1,\cdots,x_n\}.
\]
In order to obtain the relation it is necessary to attach 2-cells
\[
 \{y_1,\cdots,y_{n-1}\},
\]
where $y_j$ realizes the relator $r_j$. Let $\tilde{L}$ be the
universal covering of $L$ and $L_{\infty}$ an infinite cyclic covering
which corresponds ro ${\rm Ker}\epsilon$. According to $p=0$
(resp. $p=1$ or $p=2$), the $p$-th chain group
$C_{p}(\tilde{L},{\mathbb C})$ is a free right ${\mathbb C}[\pi_1(X_K)]$
module generated by $P_0$ (resp. $\{x_1.\cdots.x_n\}$ or
$\{y_1,\cdots,y_{n-1}\}$) and the chain complex $C_{\cdot}(L_{\infty},\rho)$ is defined to be
\[
 C_{p}(L_{\infty},\rho)=C_{p}(\tilde{L},{\mathbb C})\otimes
 _{{\mathbb C}[{\rm Ker}\epsilon]}\rho.
\] 
Thus the chain complex
\[
 C_{2}(L_{\infty},\rho)\stackrel{\partial_2}\to
 C_{1}(L_{\infty},\rho)\stackrel{\partial_1}\to C_{0}(L_{\infty},\rho),
\]
becomes
\begin{equation}
(\Lambda^{\oplus r})^{n-1} \stackrel{\partial_2}\to (\Lambda^{\oplus
r})^{n} \stackrel{\partial_1}\to \Lambda^{\oplus r},
\end{equation}
and the differentials are discribed by Fox's free differential
calculus. In fact it is known (\cite{KL}): 
\[
\partial_1=
\left(
\begin{array}{c}
\Phi(x_1-1)\\
\vdots\\
\Phi(x_n-1)
\end{array}
\right)
=
\left(
\begin{array}{c}
\rho(x_1)t-I_r\\
\vdots\\
\rho(x_n)t-I_r
\end{array}
\right),
\]
and 
\[
\partial_2=
\left(
\begin{array}{ccc}
\Phi(\frac{\partial r_1}{\partial x_1})& \cdots & \Phi(\frac{\partial
 r_1}{\partial x_n})\\
\vdots &\ddots & \vdots\\
\Phi(\frac{\partial r_{n-1}}{\partial x_1})& \cdots & \Phi(\frac{\partial
 r_{n-1}}{\partial x_n})
\end{array}
\right).
\]
Here an each entry is an element of
$M_r(\Lambda)$. $C_{p}(L_{\infty},\rho)$ is considered as a space of row vectors and
differentials act from the right. 
It is known that the determinant of a certain entry of $\partial_1$ is not
zero(\cite{Wada}). Therefore rearranging indices we may assume that
$\det(\rho(x_n)t-I_r)$ is not zero and will denote it by
$\Delta_0(t)$. Let us put
\[
 \Delta_{1}(t)={\rm det}
\left(
\begin{array}{ccc}
\Phi(\frac{\partial r_1}{\partial x_1})& \cdots & \Phi(\frac{\partial
 r_1}{\partial x_{n-1}})\\
\vdots &\ddots & \vdots\\
\Phi(\frac{\partial r_{n-1}}{\partial x_1})& \cdots & \Phi(\frac{\partial
 r_{n-1}}{\partial x_{n-1}}).
\end{array}
\right).
\]
Then {\it the twisted Alexander function} is defined to be (\cite{KL}\cite{Kitano}\cite{Wada}),
\[
 \Delta_{K,\rho}(t)=\frac{\Delta_1(t)}{\Delta_0(t)}.
\]
In the following we will assume
$\Delta_1(t)$ is not zero. This implies that after tensored with
${\mathbb C}(t)$ (30) becomes acyclic. Thus $H_{\cdot}(L_{\infty},\rho)$
are torsion $\Lambda$-modules and in particular they are finite
dimensional vector spaces. Let $\tau_i$ be the representation matrix of
$t$. Remember that $C_{\cdot}(L,\rho)$ is quasi-isomorphic to
$C_{\cdot}(X_K,\rho)$ and that
\[
 C_{\cdot}(L,\rho)=C_{\cdot}(L_{\infty},\rho)\otimes_{\Lambda} {\mathbb C}.
\]
Here ${\mathbb C}$ is regarded as a $\Lambda$-module by ${\mathbb C}\simeq
\Lambda/(t-1)$.
Thus the exact sequence of complexes
\[
 0\to C_{\cdot}(L_{\infty},\rho)\stackrel{\cdot (t-1)}\to
 C_{\cdot}(L_{\infty},\rho)\to C_{\cdot}(L,\rho) \to 0,
\]
induces
\begin{equation}
 \begin{array}{cccccccccc}
0&\to & H_2(L_{\infty},\rho)&\stackrel{\tau_2-id}\to & H_2(L_{\infty},\rho) &\to &
 H_2(X_K,\rho) & & \\
& \to & H_1(L_{\infty},\rho)&\stackrel{\tau_1-id}\to & H_1(L_{\infty},\rho) &\to &
 H_1(X_K,\rho) & & \\
& \to & H_0(L_{\infty},\rho) & \stackrel{\tau_0-id}\to &  H_0(L_{\infty},\rho) &\to & H_0(X_K,\rho) &\to & 0.
\end{array}
\end{equation}
Since $h^0(X_K,\rho)=0$,  $H_0(X_K,\rho)$ vanishes by the
universal coefficient theorem. Thus $\tau_0-id$ is an isomorphism. Using
the fact $\Delta_i(t)$ is a multiplication of the characteristic
polynomial of $\tau_i$ and a certain unit of $\Lambda$, the
following lemma is an easy consequence of (31).
\begin{lm}
The following are equivalent:
\begin{enumerate}
\item $\Delta_{K,\rho}(1)\neq 0$.
\item $h^1(X_K,\rho)=0$.
\item $h^1(X_K,\rho)=h^2(X_K,\rho)=0$.
\end{enumerate}
\end{lm}
In \cite{Sugiyama} we have proved that the fact $h^{i}(X_K,\rho)=0$ for all $i$ implies
\begin{equation}
 \tau(X_K,\rho)=|\Delta_{K,\rho}(1)|.
\end{equation}
Thus {\bf Corollary 3.4} and (32) prove the following.
\begin{thm}
Suppose $\Delta_{K,\rho}(1)\neq 0$. Then 
\[
 R_{X_K}(0,\rho)=|\Delta_{K,\rho}(1)|^2.
\]
\end{thm}

Here is an example of $\rho$ such that a special value of Ruelle
L-function at the origin can be computed explicitly. Let $\xi$ be a
complex number of modulus one. We define a morphism 
\[
 H_{1}(X_K,{\mathbb Z})\stackrel{\rho}\to U(1)
\]
to be
\[
 \rho(t)=\xi.
\]
Composing it with Hurewicz map
we obtain a unitary character
\[
 \pi_1(X_K)\stackrel{\rho}\to U(1),
\]
which yields a unitary local system of rank one on $X_K$. 
Notice that  $t$ represents a meridian of the boundary of a tubular
neighborhood of $K$. Thus if $\xi\neq 1$, $\rho$ is cuspidal. Moreover
it is known
(\cite{KL}\S 3.3):
\[
 \Delta_0(t)=1-\xi t, \quad \Delta_1(t)=A_K(\xi t),
\]
where $A_{K}(t)$ is the Alexander polynomial. Now let us choose $\xi$
so that $\xi\neq 1$ and that $A_{K}(\xi)\neq 0$. By {\bf Theorem 3.6} we have the following.
\begin{cor}
\[
 R_{X_K}(0,\rho)=\left|\frac{A_K(\xi)}{1-\xi}\right|^2.
\]
\end{cor}

\section{A geometric meaning of coefficients}
\subsection{The K-group and regulators}
\begin{fact}(\cite{Goncharov}) Let
\[
 K_{2n+1}({\mathbb C})\stackrel{r_{n+1}}\to {\mathbb R}
\]
 be the Borel regulator map. Then there is a natural element $\gamma_X$
 in $K_{2n+1}(\overline{{\mathbb Q}})\otimes{\mathbb Q}$ such that
\[
 vol(X)= r_{n+1}(\gamma_X).
\]
\end{fact}
Thus {\bf Theorem 2.2} shows if $n$ is odd the ratio of the first and
the second coefficient of the Taylor expansion at the origin is
interpretated as an evaluation of the Borel regulator against a certain element of the algebraic
K-group whereas if $n$ is even a correction from
cusps is neccessary. It seems natural to expect the first
coefficient also has such an interpretation. In fact it is true at least
for a
hyperbolic threefold, which will be explained below. Following
\cite{Dupont} we will also
explain how to construct $\gamma_X$ for a closed hyperbolic threefold.\\

Let $M$ be a smooth manifold and $P\to M$ a principal ${\mathbb C}^{\times}$-bundle with a flat
connection.  By Chern-Weil theory the image of the first Chern class $c_1(P)\in
H^{2}(M,{\mathbb Z})$ in $H^{2}(M,{\mathbb C})$ vanishes. Thus the
exact sequence
\[
 H^{1}(M,{\mathbb C}/{\mathbb Z})\stackrel{\beta}\to H^{2}(M,{\mathbb Z}) \to H^{2}(M,{\mathbb C})
\]
shows there is an element $\hat{c}_1(P)\in H^{1}(M,{\mathbb C}/{\mathbb
Z})$ which maps to $c_1(P)$ by $\beta$. Let ${\rm GL}_1({\mathbb
C})^{\delta}$ be the multiplicative group ${\mathbb C}^{\times}$ with the discrete topology and $B{\rm GL}_1({\mathbb
C})^{\delta}$ the classfying space. If we apply the previous
construction to the universal flat ${\mathbb C}^{\times}$-bundle, we
obtain an element $\hat{c}_1$ of $H^{1}(B{\rm GL}_1({\mathbb
C})^{\delta},{\mathbb C}/{\mathbb Z}).$  Since 
\[
 H^{1}(B{\rm GL}_1({\mathbb
C})^{\delta},{\mathbb C}/{\mathbb Z})\simeq H^{1}({\rm GL}_1({\mathbb
C}),{\mathbb C}/{\mathbb Z})\simeq {\rm Hom}({\mathbb C}^{\times}, {\mathbb C}/{\mathbb Z}),
\]
$\hat{c}_1$ may be regarded as a homomorphism from $K_1({\mathbb
C})\simeq {\mathbb
C}^{\times}$ to ${\mathbb C}/{\mathbb Z}$. By definition $r_1$ is
${\rm Im}\hat{c}_1$ and it is known
\[
 \hat{c}_1(g)=\frac{i}{2\pi}\log g,\quad g\in {\mathbb
C}^{\times},
\]
and thus $r_1(g)=\log |g|/2\pi$. Let us choose a unitary basis ${\bf
e}=\{{\bf e}_1,\cdots,{\bf e}_r\}$ of $\rho$. Suppose that $H^1(X,\rho)$
vanishes. Then as we have seen at the beginning of \S3.4, $C^{\cdot}(X,\rho)$ is acyclic. 
Following
\cite{MilnorI} a certain $\tau(X,\rho,{\bf e})\in K_1({\mathbb C})$ is
defined which will be referred as {\it the Milnor
element }. The Franz-Reidemeister torsion is nothing but its image in
$K_1({\mathbb C})/{\rm U}(1)\simeq {\mathbb R}$, i.e. its modulus. 
Thus 
\begin{eqnarray*}
\log R_X(0,\rho) &=& 2\log \tau(X,\rho)\\
&=&2 \log |\tau(X,\rho,{\bf e})|\\
&=& 4\pi r_1(\tau(X,\rho,{\bf e})). 
\end{eqnarray*}
and we have found that $\log R_X(0,\rho)$ is interpretated as a period
of the Milnor elememt by the first Borel regulator.\\

Taking account of an exceptional isomorphism
\[
 {\rm Spin}(3,1)\simeq {\rm SL}_2({\mathbb C}), 
\]
let us apply the previous construction to the universal flat bundle ${\rm
SL}_2({\mathbb C})$-bundle on $B{\rm
SL}_2({\mathbb C})^{\delta}$. Then we will obtain the Chern-Simon class
\begin{eqnarray*}
 \hat{c}_2&\in& H^3(B{\rm
SL}_2({\mathbb C})^{\delta},{\mathbb C}/{\mathbb Z})\simeq H^3({\rm
SL}_2({\mathbb C}),{\mathbb C}/{\mathbb Z})\\
&\simeq & {\rm Hom}(H_3({\rm SL}_2({\mathbb C}),{\mathbb Z}),{\mathbb
 C}/{\mathbb Z})
\end{eqnarray*} 
which is a lift of the second Chern class $c_2\in H^4(B{\rm
SL}_2({\mathbb C})^{\delta},{\mathbb Z})$. In particular its imaginary
part ${\rm Im}\hat{c}_2$ yields a ${\mathbb Q}$-linear map:
\[
 H_3({\rm SL}_2({\mathbb C}),{\mathbb Q})\stackrel{{\rm Im}\hat{c}_2}\to
 {\mathbb R},
\]
and it is known (\cite{Dupont})
\[
 {\rm Im}\hat{c}_2(g_1|g_2|g_3)=\frac{1}{4\pi^2}vol(T(\infty, g_1\infty,
 g_1g_2\infty, g_1g_2g_3\infty)),\quad \,g_i\in {\rm SL}_2({\mathbb C}).
\]
Here $T(z_1,z_2,z_3,z_4)$ is a tetrahedron in ${\mathbb H}^3$ whose vertices are
$\{z_1,z_2,z_3,z_4\}$ and edges are geodesics. Remember that the
volume of an ideal tetrahedron is computed by {\it the Bloch-Wigner
function}:
\[
 D(z)={\rm arg}(1-z){\rm arg}(z)+{\rm Im}{\rm Li}_2(z),
\]
where ${\rm Li}_2(z)$ is {\it the dilogarithm}:
\[
 {\rm Li}_2(z)=\sum_{n=1}^{\infty}\frac{z^n}{n^2}.
\]
More precisely it is known (\cite{Dupont})
\begin{eqnarray*}
vol(T(\infty, g_1\infty,g_1g_2\infty,g_1g_2g_3\infty))&=& vol(T(\infty,
 0,1,z(g_1,g_2,g_3)))\\
&=& D(z(g_1,g_2,g_3)),
\end{eqnarray*}
where $z(g_1,g_2,g_3)$ is the cross ratio of $\{\infty,
g_1\infty,g_1g_2\infty,g_1g_2g_3\infty\}$. Suppose that $X$ is
closed. Then since $H_3(B{\rm
PSL}_2({\mathbb C})^{\delta},{\mathbb Q})$ is isomorphic to $H_3(B{\rm
SL}_2({\mathbb C})^{\delta},{\mathbb Q})$ the natural inclusion $\Gamma \hookrightarrow {\rm SO}(3,1)\simeq
{\rm PSL}_2({\mathbb C})$ induces 
\[
 H_3(X,{\mathbb Q})\simeq H_3(B\Gamma,{\mathbb Q}) \to H_3(B{\rm
SL}_2({\mathbb C})^{\delta},{\mathbb Q})\simeq H_3({\rm
SL}_2({\mathbb C}),{\mathbb Q}).
\]
Here recall that $H_3({\rm
SL}_2({\mathbb C}),{\mathbb Z})$ is a direct summand of the Quillen's
K-group $K_3({\mathbb C})$. In fact (\cite{Dupont},(9.6)),
\[
 K_3({\mathbb C})\simeq H_3({\rm
SL}_2({\mathbb C}),{\mathbb Z})\oplus K_3^{M}({\mathbb C}),
\]
where $K_3^{M}({\mathbb C})$ is the Milnor K-group. Thus the fundamental
    class $[X]$ of $X$ defines an element $\gamma_X$ in $K_3({\mathbb
    C})\otimes{\mathbb Q}$. According to Weil's rigidity $\Gamma$ is conjugate a subgroup of ${\rm
    PSL}_2(\overline{\mathbb Q})$ and   
    $\gamma_X$ is contained in $K_3(\overline{\mathbb
    Q})\otimes{\mathbb Q}$. We define {\it the second Borel regulator}
\[
 K_3({\mathbb
    C})\otimes{\mathbb Q}\stackrel{r_2}\to {\mathbb R}
\]
to be a composition of $4\pi^2{\rm Im}\hat{c}_2$ and the natural
projection:
\[
 K_3({\mathbb
    C})\otimes{\mathbb Q}\to H_3({\rm
SL}_2({\mathbb C}),{\mathbb Q}).
\]
Then {\bf Corollary 2.3} yields 
\[
 \frac{d}{dz}\log R_X(z,\rho)|_{z=0}=-\frac{3r}{\pi}\cdot r_2(\gamma_X).
\]
If $X$ is not closed, $H_3(X)$ vanishes. So we have to use the relative
homology group to define $\gamma_X$. (See \cite{Goncharov} for details.)
Thus we have found that the leading and the second coefficient of Taylor
    expansion of $R_X(z,\rho)$ at the origin are expressed by the logarithm
    and the dilogarithm, respectively.

\subsection{The $L^2$-torsion}
 The constant term of the logarithmic derivative of Ruelle L-function at
the origin is also related to $L^2$-{\it analytic torsion}
(\cite{Atiyah1976}\cite{Mathai}). Following \cite{Atiyah1976}, we
remember the von Neumann trace. Let $\omega_p$ the action of $\Gamma$ on
$L^2({\mathbb H}^d, \Omega^p)$. Then $L^2({\mathbb H}^d,
\Omega^p\otimes{\mathbb C}^r)$ is a $\Gamma$-module by
$\omega_p\otimes\rho$ and there is an isomorphism of $\Gamma$-modules:
\begin{equation}
 L^2({\mathbb H}^d,
\Omega^p\otimes{\mathbb C}^r)\simeq L^2(\Gamma)\otimes L^2({\mathbb H}^d, \Omega^p).
\end{equation}
Here $\Gamma$ acts on $L^2(\Gamma)$ by the left regular representaion
and we regard $L^2({\mathbb H}^d,
\Omega^p)$ is the trivial module. Since Hodge Laplacian
$\Delta^p$ on $L^2({\mathbb H}^d,
\Omega^p\otimes{\mathbb C}^r)$ commutes with $\Gamma$ so does
$e^{-t\Delta^p}$. Let $U$ be the fundamental domain of $\Gamma$ and
$\psi$ its characteristic function. Using (14) von Neumann trace
of $e^{-t\Delta^p}$ is given by (\cite{Atiyah1976}) 
\[
 \tau(e^{-t\Delta^p})={\rm Tr}(\psi\cdot e^{-t\Delta^p}\cdot\psi),
\]
which is equal to $I_p(t)$. Let us put
\[
 \zeta_2(s)=\sum_p(-1)^pp\cdot \zeta_2^{(p)}(s),\quad \zeta_2^{(p)}(s)=\frac{1}{\Gamma(s)}M(\tau(e^{-t\Delta^p}))(s).
\]

Then {\it the analytic $L^2$-torsion} of $(X,\rho)$ is defined to be
\[
 \tau_{an}^{(2)}(X,\rho)=\exp(-\frac{1}{2}\zeta^{\prime}_2(0)).
\]
\cite{Mathai}{\bf Lemma 6.4} (see also {\bf Appendix}, {\bf Lemma 5.1}
below) and the computation in  \S2.2 show
that 
$M(I_p)(0)$ is a rational multiple of $vol(X)/\pi$ and that 
\[
 \frac{d}{ds}\zeta_2^{(p)}(s)|_{s=0}=M(I_p)(0).
\]

Now {\bf Theorem 2.2} implies
\begin{thm}
Let $h$ be the order of $R_X(\rho,z)$ at the origin. Then there is a
 rational number $\alpha$ such that 
\[
 \lim_{z\to
 0}\{\frac{d}{dz}\log R_X(\rho,z)-\frac{h}{z}\}-2\sum_{j=0}^n(-1)^n\delta(X,\rho)=\alpha\cdot
 \log \tau^{(2)}_{an}(X,\rho).
\]
\end{thm}
For example suppose $d=3$.  Since (\cite{Mathai} {\bf Corollary 6.7})
\[
 \log \tau^{(2)}_{an}(X,\rho)=M(I_0)(0)=\frac{r}{6\pi}vol(X),
\]
we obtain
\[
 \lim_{z\to
 0}\{\frac{d}{dz}\log R_X(\rho,z)-\frac{2h^1(X,\rho)}{z}\}=-18\cdot
 \log \tau^{(2)}_{an}(X,\rho).
\]
\section{Appendix}
We will show {\bf Fact 3.1} under the assumption that $\rho$ is cuspidal and $d=3$. In the following the Mellin transform of a
function $f$ will be denoted by $M(f)$:
\[
 M(f)(s)=\int^{\infty}_0f(t)t^{s-1}dt.
\]
\begin{lm} Let $t$ be a positive number.
\begin{enumerate}
\item 
\[
 \int^{\infty}_{-\infty}e^{-t\lambda^2}d\lambda=\sqrt{\pi}t^{-\frac{1}{2}}.
\]
\item For a positive integer $k$,
\[
 \int^{\infty}_{-\infty}e^{-t\lambda^2}\lambda^{2k}d\lambda=\frac{\sqrt{\pi}(2k-1)!!}{2^{k}}t^{-\frac{1}{2}-k}.
\]
\item Let $c$ be a positive number and $P$ an even polynomial. Then the
      Mellin transform of 
      $\int^{\infty}_{-\infty}e^{-t(\lambda^2+c^2)}P(\lambda)
      d\lambda$ is meromorphically continued to ${\mathbb C}$. It is regular at $s=0$ and 
\[
 M(\int^{\infty}_{-\infty}e^{-t(\lambda^2+c^2)}P(\lambda)d\lambda)(0)=-2\pi \int^{c}_0P(iy)dy.
\]
\end{enumerate}
\end{lm}
{\bf Proof.} See \cite{Fried}{\bf Lemma 3}.
\begin{flushright}
$\Box$
\end{flushright}
We put
\[
 \theta_p(t)={\rm Tr}[e^{-\Delta^p_X}]-h^p(X,\rho).
\]
Since $H^0(X,\rho)$ vanishes Selberg trace formula shows
\begin{equation}
 \theta_0(t)=\delta_0(t)=h_0(t)+e_0(t),
\end{equation}

where 
\begin{eqnarray*}
 e_0(t)&=& i_0(t)+u_0(t)\\
&=& \frac{r\cdot
 vol(X)}{4\pi^2}\int^{\infty}_{-\infty}e^{-t(\lambda^2+1)}\lambda^2d\lambda+\frac{\delta(X,\rho)}{2\pi}\int^{\infty}_{-\infty}e^{-t(\lambda^2+1)}d\lambda\\
&=& \frac{r\cdot vol(X)}{8\pi\sqrt{\pi}}e^{-t}t^{-\frac{3}{2}}+\frac{\delta(X,\rho)}{2\sqrt{\pi}}e^{-t}t^{-\frac{1}{2}}.
\end{eqnarray*}
By {\bf Lemma 5.1}, $M(e_0)(s)$ is regular at the origin and 
\[
 M(e_0)(0)=\frac{r}{6\pi}vol(X)-\delta(X,\rho).
\]
If $t>0$ is sufficiently small, 
\[
 h_0(t)\sim \frac{a}{\sqrt{4\pi t}}e^{-\frac{c_X^2}{4t}},
\] 
where $c_X$ is the length of minimal closed geodesic.  Since
$\theta_0(t)$ exponentially decays as $a\to \infty$ so does $h_0(t)$. Therefore $M(h_0)(s)$ is an entire function and $M(\theta_0)(s)$ is
a meromorphic function on the whole plane regular at the origin. Notice
that $\Gamma(s)$ has simple pole with residue $1$ at $s=0$.
\begin{prop}
$\zeta_X^{(0)}(s,\rho)=M(\theta_0)(s)/\Gamma(s)$ satisfies the following properties.
\begin{enumerate}
\item It is a meromorphic function on ${\mathbb C}$ and vanishes at the
      origin. 
\item 
\[
 \frac{d}{ds}\zeta_X^{(0)}(s,\rho)|_{s=0}=M(\theta_0)(0)=M(h_0)(0)+\frac{r}{6\pi}vol(X)-\delta(X,\rho).
\]
\end{enumerate}
\end{prop}
By definition the functional determinant is 
\[
 \det \Delta_X^{p}=\exp(-\frac{d}{ds}\zeta_X^{(p)}(s,\rho)|_{s=0}).
\]
Hence
\begin{equation}
-\log\det\Delta_X^{0}=M(\theta_0)(0)=M(h_0)(0)+\frac{r}{6\pi}vol(X)-\delta(X,\rho).
\end{equation}
We put
\[
 \eta_1(t)=h_1(t)+e_1(t)-h^1(X, \rho),
\]
where 
\begin{eqnarray*}
 e_1(t)&=& i_1(t)+u_1(t)\\
&=& \frac{r\cdot
 vol(X)}{\pi^2}\int^{\infty}_{-\infty}e^{-t\lambda^2}(\lambda^2+1)d\lambda+\frac{\delta(X,\rho)}{\pi}\int^{\infty}_{-\infty}e^{-t\lambda^2}d\lambda\\
&=& \frac{r\cdot vol(X)}{2\pi\sqrt{\pi}}(2t^{-\frac{1}{2}}+t^{-\frac{3}{2}})+\frac{\delta(X,\rho)}{\sqrt{\pi}}t^{-\frac{1}{2}}.
\end{eqnarray*}
Then
\[
 \theta_1(t)=\eta_1(t)+\theta_0(t).
\]
In order to investigate the Mellin transform of $\eta_1$ we consider
\begin{eqnarray}
\mu_1(t)&=& \eta_1(t) +(h^1(X, \rho)-e_1(t))\cdot \chi_{(0,1]}\\
&=& h_1(t) -(h^1(X, \rho)-e_1(t))\cdot \chi_{(1,\infty)},
 \end{eqnarray}
where $\chi$ is a characteristic function. For a sufficiently small
positive number $t$, (37) shows
\[
 \mu_1(t) \sim h_1(t) \sim \frac{a}{\sqrt{4\pi t}}e^{-\frac{c_X^2}{4t}}
\]
and (36) implies for sufficiently large $t$
\[
 \mu_1(t) \sim \eta_1(t) \sim e^{-\gamma t},\quad \gamma > 0.
\]
Thus $M(\mu_1)(s)$ is an entire function.  If ${\rm
Re} s>\frac{3}{2}$,
\begin{eqnarray*}
M(\eta_1)(s) &=& M(\mu_1)(s)+\int^{1}_0(e_1(t)-h^1(X,\rho))t^{s-1}dt\\
&=& M(\mu_1)(s)-\frac{h^1(X,\rho)}{s}+(\frac{r\cdot vol(X)}{\pi\sqrt{\pi}}+\frac{\delta(X,\rho)}{\sqrt{\pi}})\frac{1}{s-1/2}+\frac{r\cdot vol(X)}{2\pi\sqrt{\pi}}\frac{1}{s-3/2}.
\end{eqnarray*}
and if ${\rm Re} s<\frac{1}{2}$,
\begin{eqnarray*}
M(h_1)(s) &=& M(\mu_1)(s)-\int^{\infty}_1(e_1(t)-h^1(X,\rho))t^{s-1}dt\\
&=& M(\mu_1)(s)-\frac{h^1(X,\rho)}{s}+(\frac{r\cdot vol(X)}{\pi\sqrt{\pi}}+\frac{\delta(X,\rho)}{\sqrt{\pi}})\frac{1}{s-1/2}+\frac{r\cdot vol(X)}{2\pi\sqrt{\pi}}\frac{1}{s-3/2},
\end{eqnarray*}
we see that both $M(\eta_1)(s)$ and $M(h_1)(s)$ are meromorphically
continued to the whole plane as the same function.
Moreover we find that $M(\eta_1)(s)+h^1(X,\rho)/s$ is regular at the
origin. Together with {\bf Proposition 5.1} this shows
\begin{prop}
$\zeta_X^{(1)}(s,\rho)=M(\theta_1)(s)/\Gamma(s)$ satisfies the following properties.
\begin{enumerate}
\item It is a meromorphic function on ${\mathbb C}$ and is regular at the origin. Moreover
\[
 \zeta_X^{(1)}(0,\rho)=-h^1(X,\rho).
\] 
\item 
\[
 \frac{d}{ds}\zeta_X^{(1)}(s,\rho)|_{s=0}=M(\theta_0)(0)+\lim_{s\to
  0}\{\Gamma(s)h^1(X,\rho)+M(\eta_1)(s)\}.
\]
\end{enumerate}
\end{prop}
By {\bf Proposition 2.1}, {\bf Corollary 2.1}, {\bf Proposition 2.4} and
(34) we obtain
\begin{eqnarray*}
 s_0(z+1)&=& L(e^th_0)(z)\\
&=& L(e^t\delta_0)(z)-L(e^ti_0)(z)-L(e^tu_0)(z)\\
&=& L(e^t\delta_0)(z)+\frac{r\cdot vol(X)}{2\pi}z^2-\delta(X,\rho).
\end{eqnarray*}
By {\bf Lemma 2.6} $L(e^t\delta_0)(z)$ is an odd function and thus
\[
 s_0(1-z)+s_0(1+z)=\frac{r}{\pi}vol(X)z^2-2\delta(X,\rho).
\]
Since $s_j(z)$ is the loratithmic derivative of $S_j(z)$, (35) yields
\begin{eqnarray*}
 \log S_0(2)-\log S_0(0) &=& \int^1_0(s_0(1+z)+s_0(1-z))dz\\
&=& \frac{r\cdot vol(X)}{3\pi}-2\delta(X,\rho)\\
&=& -2\log\det\Delta^0_X-2M(h_0)(0).
\end{eqnarray*}
The equation (\cite{Fried}, pp535 (13)):
\[
 M(h_0)(0)=-\log S_0(2).
\]
and (35) implies
\begin{prop}
\[
 \log (S_0(0)S_0(2))=2\log\det\Delta^0_X =-2M(\theta_0)(0).
\]
\end{prop}

Now we will compute the Ray-Singer torsion. By {\bf Proposition 5.1},
{\bf Proposition 5.2} and {\bf Proposition 5.3} 
\begin{eqnarray*}
\zeta_X^{\prime}(0,\rho)&=&
 \frac{d}{ds}\zeta_X^{(1)}(s,\rho)|_{s=0}-3\frac{d}{ds}\zeta_X^{(0)}(s,\rho)|_{s=0}\\
&=& \lim_{s\to 0}\{\Gamma(s)h^1(X,\rho)+M(\eta_1)(s)\}+\log (S_0(0)S_0(2)).
\end{eqnarray*}
The arguments of pp.536 of \cite{Fried}(especially (25)) shows that the leading term
of the Taylor expansion
of $S_1(z+1)$ at the origin is $\delta z^{2h^1(X,\rho)}$. Here $\delta$
is given by
\[
 -\log \delta =\lim_{s\to 0}\{\Gamma(s)h^1(X,\rho)+M(\eta_1)(s)\}.
\]
Thus 
\[
 \lim_{s\to 0}\{\Gamma(s)h^1(X,\rho)+M(\eta_1)(s)\}=-\log(\lim_{z\to
 0}z^{-2h^1(X,\rho)}S_1(z+1)), 
\]
and {\bf Fact 2.1} shows
\begin{thm}
\[
 \lim_{z\to 0}z^{-2h^1(X,\rho)}R_X(z,\rho)=\exp(-\zeta_X^{\prime}(0,\rho)).
\]
\end{thm}


\begin{thebibliography}{10}

\bibitem{Atiyah1976}
M.~F. Atiyah.
\newblock Elliptic operators, discrete groups and von {N}eumann algebras.
\newblock {\em Ast\'{e}risque}, 32-33:43--72, 1976.

\bibitem{Bismut-Zhang}
J.~M. Bismut and W.~Zhang.
\newblock {\em An extension of a theorem by {C}heeger and {M}\"{u}ller}, volume
  205 of {\em Ast\'{e}risque}.
\newblock S.M.F., 1992.

\bibitem{Dai-Fang}
X.~Dai and H.~Fang.
\newblock Analytic torsion and {R}-torsion for manifolds with boundary.
\newblock {\em Asian J. Math.}, 4(3):695--714, September 2000.

\bibitem{Dupont}
J.~L. Dupont.
\newblock {\em Scissors Congruences, Group Homology and Characteristic
  Classes}.
\newblock World Scientific, 2001.

\bibitem{Fried}
D.~Fried.
\newblock Analytic torsion and closed geodesics on hyperbolic manifolds.
\newblock {\em Inventiones Math.}, 84:523--540, 1986.

\bibitem{Goncharov}
A.~B. Goncharov.
\newblock Volume of hyperbolic manifolds and mixed {T}ate motives.
\newblock {\em J. Amer. Math.Soc.}, 12:569--618, 1999.

\bibitem{KL}
P.~Kirk and C.~Livingston.
\newblock Twisted {A}lexander invarinants, {R}eidemeister torsion, and
  {C}asson-{G}ordon invariants.
\newblock {\em Topology}, 38(3):635--661, 1999.

\bibitem{Kitano}
T.~Kitano.
\newblock Twisted {A}lexander polynomial and {R}eidemeister torsion.
\newblock {\em Pacific Jour. Math.}, 174(2):431--442, 1996.

\bibitem{Mathai}
V.~Mathai.
\newblock ${L}^2$-analytic torsion.
\newblock {\em J. Funct. Anal.}, 107:369--386, 1992.

\bibitem{Miatello}
R.~J. Miatello.
\newblock On the {P}lancherel measure for linear {L}ie groups of rank one.
\newblock {\em Manuscripta Math.}, 29(2-4):249--276, 1979.

\bibitem{MilnorW}
J.~Milnor.
\newblock Whitehead torsion.
\newblock {\em Bull. Amer. Math. Soc.}, 72:358--426, 1966.

\bibitem{MilnorI}
J.~Milnor.
\newblock Infinite cyclic coverings.
\newblock In J.~G. Hocking, editor, {\em Conference on the Topology of
  Manifolds}, pages 115--133. PWS Publishing Company, 1968.

\bibitem{Osborne-Warner}
M.~S. Osborne and G.~Warner.
\newblock Multiplicities of the integral discrete series : the case of a
  nonuniform lattice in an {R}-rank one semisimple group.
\newblock {\em J. Funct. Anal.}, 30(3):287--310, 1978.

\bibitem{Park2007}
J.~Park.
\newblock Analytic torsion and closed geodesics for hyperbolic manifolds with
  cusps.
\newblock Preprint, February 2008.

\bibitem{Ray-Singer}
D.~B. Ray and I.~M. Singer.
\newblock R-torsion and {L}aplacian on {R}iemannian manifolds.
\newblock {\em Advances in Math.}, 7:145--210, 1971.

\bibitem{Reed-Simon}
M.~Reed and B.~Simon.
\newblock {\em Methods of Modern Mathematical Physics}, volume~4.
\newblock Academic Press, 1978.
\newblock Analysis of Operators.

\bibitem{Sarnak-Wakayama}
P.~Sarnak and M.~Wakayama.
\newblock Equidistribution of holonomy about closed geodesics.
\newblock {\em Duke J. Math.}, 100(1-57), 1999.

\bibitem{Sugiyama}
K.~Sugiyama.
\newblock An analog of the {I}wasawa conjecture for a compact hyperbolic
  threefold.
\newblock {\em J. Reine Angew. Math.}, 613:35--50, 2007.

\bibitem{Wada}
M.~Wada.
\newblock Twisted {A}lexander polynomials for finitely presented groups.
\newblock {\em Topology}, 33(2):241--256, 1994.

\end{thebibliography}

\vspace{10mm}
\begin{flushright}
Address : Department of Mathematics and Informatics\\
Faculty of Science\\
Chiba University\\
1-33 Yayoi-cho Inage-ku\\
Chiba 263-8522, Japan \\
e-mail address : sugiyama@math.s.chiba-u.ac.jp
\end{flushright}

\end{document}